\documentclass[a4paper,12pt]{article}
\usepackage[english, french]{babel}
\usepackage{amsfonts, amsmath,amscd,amssymb} 
\usepackage{latexsym}
\usepackage[latin1]{inputenc}
\usepackage{url}
\makeatletter
\@addtoreset{equation}{section}

\makeatother
 \usepackage[bookmarks]{hyperref}
 \usepackage{tikz}

\newcommand{\C}{\mathbb{C}}
\newcommand{\R}{\mathbb{R}}
\newcommand{\N}{\mathbb{N}}

\newcommand{\im}{\mathrm{Im}\,}
\newcommand{\re}{\mathrm{Re}\,}
\newcommand{\MO}{\mathcal{O}}

\begin{document}
\title{Estimation de résolvante et construction de quasimode près du bord
du pseudospectre}
\author{William Bordeaux Montrieux \\ 
\small Department of Mathematics\\ 
\small University of Toronto\\
\small Toronto, Ontario, Canada, M5S 2E4 \\
\small william.montrieux@utoronto.ca}
\date{}
\maketitle

\newtheorem{theo}{Th\'{e}orème}[section]
\newtheorem{corol}[theo]{Corollaire}
\newtheorem{prop}[theo]{Proposition}
\newtheorem{lemme}[theo]{Lemme}
\newtheorem{definition}[theo]{D\'{e}finition}
\newtheorem{hypothese}[theo]{Hypothèse}
\newtheorem{remarque}[theo]{Remarque}
\newtheorem{conjecture}[theo]{Conjecture}
\newtheorem{exemple}[theo]{Exemple}

\begin{abstract}
Nous considérons un opérateur $h$-pseudodifférentiel non-autoadjoint
dans la limite semiclassique.  $p$ désigne le symbole principal.
Nous savons que la résolvante
existe à l'intérieur de l'image de $p$ jusqu'à une distance 
$\MO((h\ln\frac{1}{h})^{\frac{k}{k+1}}),$ de certains points du bord, où
$k\in\{2,4,\ldots\}.$ 
Dans ce travail,
nous précisons les estimations de résolvantes qu'ont obtenues
différents auteurs dans le cas $k=2,$ et en dimension $1.$
Pour la preuve, il s'agit de construire, via un scaling, 
des quasimodes pour des valeurs
du paramètre spectral très proches du bord de l'image de $p.$
\end{abstract}
\selectlanguage{english}   
\begin{abstract} 
We consider a non-self-adjoint pseudodifferential 
operator in the semi-classical limit
$(h\to0)$.
The principal symbol is given by $p.$
We know that the resolvent $(z-P)^{-1}$ exists inside the range up to a distance
$\MO((h\ln\frac{1}{h})^{\frac{k}{k+1}})$
from certain boundary points, where $k\in\{2,4,\ldots\}.$
In this work, we improve the resolvent estimates
given by different authors in the case $k=2,$ and in 
dimension one. For the proof, we will construct  quasimodes
by a scaling for  $z$ very close to the boundary. 
 \end{abstract}
\selectlanguage{french}  

\tableofcontents
\section{Introduction}
Dans \cite{DSZ}, N.~Dencker, J.~Sj\"ostrand et M.~Zworski  obtiennent
une estimation de la résolvante pour certains points du bord de l'image 
du symbole principal. Dans le cas d'un point $z_{0}$ de type fini d'ordre $k,$ 
ils obtiennent que la résolvante a une croissance en $h^{-\frac{k}{k+1}},$ 
dans des disques de rayon $\MO(h^{\frac{k}{k+1}})$ centrés en $z_{0}.$
Dans \cite{SJ}, J.~Sj\"ostrand montre que la résolvante peut s'étendre à des disques
de rayon $\MO((h\ln\frac{1}{h})^{\frac{k}{k+1}}),$ et donne une majoration pour la résolvante.
Une majoration similaire avait été obtenue auparavant 
par l'auteur pour un opérateur modèle $hD_{x}+g(x),$
$x\in S^{1},$ avec $k=2,$ voir \cite{WB} Chapitre 4. 
On citera aussi
J. Martinet \cite{JM} 
qui obtient un encadrement de la norme de la résolvante pour l'opérateur d'Airy 
complexe. La preuve est basée sur une analyse directe du semigroupe 
après conjugaison par la transformée de  Fourier.
Dans ce travail, il s'agit d'améliorer et généraliser les estimations de la résolvante obtenues 
par ces différents auteurs. 
On se limite au cas d'un point d'ordre $2,$ et à la dimension 1.
Les estimations obtenues nous permettent de préciser  l'estimation  
de résolvante obtenue par Martinet \cite{JM} pour l'opérateur d'Airy complexe 
sur la droite 
réelle $(D_{x})^{2}+ix,$ $x\in \R,$ et d'obtenir le comportement à l'infini
des lignes de niveaux de la résolvante de l'oscillateur harmonique non-autoadjoint.
Ce dernier résultat vient compléter ceux publiés par
L.~Boulton \cite{LB}, K.~Pravda-Starov \cite{KPS}. 
Mentionnons aussi le travail \cite{DK} de Davies et Kuijlaars qui obtiennent 
le premier terme dans l'asymptotique du logarithme de la norme des projecteurs 
spectraux pour les grandes valeurs propres pour l'oscillateur harmonique non-autoadjoint.
\paragraph{Cas modèle sur $\R$ à paramètre.}
Nous considérons l'opérateur modèle non-autoadjoint dans $L^{2}(\R)$
\[P_{\alpha}=hD_{x}+g(x,\alpha)+h^{2}g_{\ast}(x,\alpha;h),\quad h\in (0,1], \quad
D_{x}=\frac{1}{i}\frac{\partial}{\partial x},\]
$g$ et $g_{\ast}$ dépendent de manière $C^{\infty}$ de $x$ et d'un paramètre $\alpha$
strictement positif et borné.
On pose $\widetilde{g}(x,\alpha;h):=g(x,\alpha)+h^{2}g_{\ast}(x,\alpha;h)$ et 
on demande qu'il existe  $C_{0}>0$ tel que 
\[ 
\re g=\MO(x^{2}+1),\]
(cette hypothèse peut facilement être affaiblie).
Pour $g_{\ast}$ nous avons le dévelop\-pement asymptotique 
en puissances de $h$ suivant
\begin{align}
g_{\ast}(x,\alpha;h)\sim g^{0}_{\ast}(x,\alpha)+hg^{1}_{\ast}(x,\alpha)+\ldots
\mbox{ dans } C^{\infty}_{0}(\R\times[0,1]).
\end{align}
L'opérateur $P_{\alpha}$ est 
muni du domaine $H_{sc}(\R):=\{u\in L^{2}(\R)|\; \|x^{2}u\|+\|hD_{x}u\|<+\infty\}.$
On note $p_{\alpha}$ le symbole principale semiclassique de $P_{\alpha},$ soit 
$p_{\alpha}=\xi+g(x,\alpha).$ Rappelons que la définition du crochet de Poisson est
\begin{equation}
\{a,b\}(x,\xi):=(a'_{\xi}b'_{x}-a'_{x}b'_{\xi})(x,\xi)=H_{a}b,
\end{equation}
pour deux fonctions $a(x,\xi),$ $b(x,\xi)$ de classe  $C^{1}(\R^{2}).$ 
Ici $H_{a}=(a_{\xi}\partial_{x}-a_{x}\partial_{\xi})$ désigne le champ Hamiltonien. 
\begin{prop}\label{int1}
Soient $\widetilde{g}$ comme ci-dessus avec 
$g(x,\alpha)\in C^{\infty}(\R\times]0,C[),$ et pour tout $h$ fixé
$g_{\ast}(x,\alpha;h) \in C^{\infty}(\R\times[0,1/C])$
pour un $C>0$ donné.
Nous rajoutons les hypothèses 
suivantes:\\ 
pour $x\ne 0,$  $g$ satisfait $\im g(x,0)>0,$ et 
\begin{align}\label{int.4}
g(0,0)&=0,\; g'_{\alpha}(0,0)=-i,\\
g'_{x}(0,0)&=0,\; \im g''_{xx}(0,0)>0.\label{int.5}
\end{align}
L'opérateur $P_{\alpha}$ est muni du domaine $H_{sc}(\R).$
Dans ces conditions,
la résolvante est définie pour tout $\alpha$
\begin{align}
(hD_{x}+\widetilde{g}(x,\alpha;h))^{-1}: L^{2}(\R)\to L^{2}(\R),
\end{align}
et vérifie pour $h\ll \alpha^{3/2}$ et $\alpha$ assez petit,  
\begin{equation}\label{int.6}
\|P_{\alpha}^{-1}\|\sim  \frac{\sqrt{\pi}\exp(\frac{1}{h}\im \ell_{0}(\alpha))}{h^{1/2}
(\frac{1}{2i}\{p_{\alpha},\overline{p_{\alpha}}\}(\rho_{+}))^{\frac{1}{4}}
(\frac{1}{2i}\{\overline{p_{\alpha}},p_{\alpha}\}(\rho_{-}))^{\frac{1}{4}}}
(1+\mathcal{O}( \tilde{h}))
+\mathcal{O}(\frac{1}{\sqrt{h}\,\alpha^{1/4}}),
\end{equation} 
où $\tilde{h}=h/\alpha^{3/2},$ 
et $\rho_{\pm}(\alpha):=(x_{\pm}(\alpha),\xi_{\pm}(\alpha))$
sont les solutions de $p_{\alpha}(\rho_{\pm})=\xi_{\pm}+g(x_{\pm},\alpha)=0$
avec $\mp \im'(x_{\pm},\alpha)>0$
 (impliquant
$\im g(x_{\pm},\alpha)=0).$
$\ell_{0}$ est donné par 
\begin{equation}
\ell_{0}:=-\im \int_{x_{+}}^{x_{-}}g(y,\alpha)dy\asymp \alpha^{3/2}.
\end{equation}
Nous avons pour $\alpha$ assez petit
\begin{equation}
x_{\pm}(\alpha)=\mp \alpha^{\frac{1}{2}}\left(
\frac{-2}{\{p_{0},\frac{1}{2i}\{p_{0},\overline{p_{0}}\}\}(0,0)}
\right)^{1/2}
+\MO(\alpha),
\end{equation}
impliquant
\[|\frac{1}{2}\{p_{\alpha},\overline{p_{\alpha}}\}(\rho_{\pm})|
\asymp \alpha^{\frac{1}{2}}.\]
\end{prop}
Voir la section 2 pour la preuve.
\begin{remarque}
Le dénominateur du premier de \ref{int.6} est de l'ordre de $\sqrt{h}\alpha^{1/4}$ et ce
qui le rend dominant 
dès que $\alpha\gg h^{2/3}.$
\end{remarque}
\begin{remarque}\label{re3}
La proposition reste valable pour $\alpha\in[0,1]$ si on rajoute l'hypothèse 
(automatiquement vérifiée pour $\alpha$ petit) que l'équation $\im g(x,\alpha)=0$
admet exactement deux solution $x=x_{+}(\alpha)$, $x_{-}(\alpha)$ et de plus
$\mp\im g'_{x}(x_{\pm}(\alpha),\alpha)$ $>0.$
\end{remarque}
Nous considérons l'opérateur d'Airy complexe sur la droite réelle
$\mathcal{A}=D_{x}^{2}+ix,$ étudié par J. Martinet \cite{JM}, et muni du domaine
$D(\mathcal{A})=\{u\in H^{2}(\R)|\; xu\in L^{2}(\R)\}.$
Nous savons que le spectre est vide, que la norme de la résolvante 
$(\mathcal{A}-z)^{-1}$
ne dépend
que de $\re z$ (voir Helffer \cite{He}), et qu'en conjuguant par une transformée de Fourier, 
on peut se ramener à l'étude de l'opérateur
\[D_{x}-ix^{2}+iz.\]
Pour $z$ dans la zone d'intérêt $\re z>0,$
le changement de variable $x=(\re z)^{\frac{1}{2}} y,$ 
donne l'opérateur 
$(\re z)\mathcal{Q},$ où $\mathcal{Q}$ est donné par 
\[\frac{1}{(\re z)^{\frac{3}{2}}}D_{x}-ix^{2}+i-\frac{\im z}{\re z}.\]
De la proposition au-dessus et la remarque \ref{re3},
avec $h=\frac{1}{(\re z)^{\frac{3}{2}}}$ et $\alpha=1,$ on déduit 
le corollaire suivant:
\begin{corol}\label{int2}
Pour $\re z>0$ assez grand, la résolvante de $\mathcal{A}$ satisfait
\begin{align*}
\|((D_{x})^{2}+ix-z)^{-1}\|\sim &\sqrt\frac{\pi}{2}(\re z)^{-\frac{1}{4}}(1+\MO(\frac{1}{(\re z)^{3/2}}))
 \exp(\frac{4}{3}(\re z)^{\frac{3}{2}})\\
 &+\MO(\frac{1}{(\re z)^{1/4}}).
\end{align*}
\end{corol}   
Ceci précise l'estimation qu'a obtenue J. Martinet dans \cite{JM}.
La preuve sera donnée à la fin de sa sous-section 2.2.
\paragraph{Cas général en dimension 1.} 
Soit $p\in S(\R^{2},m)$ indépendant de $h,$ où $m$ désigne une fonction d'ordre
au sens que
\begin{align*}
\exists C_{0}\ge 1,\,N_{0}>0
\mbox{ tels que } m(\rho)\le C_{0}\langle \rho-\mu\rangle^{N_{0}}m(\mu),
\forall\rho,\mu\in \R^{2},\\
\mbox{ avec }
\langle \rho-\mu\rangle=\sqrt{1+|\rho-\mu|^{2}}.
\end{align*}

L'espace de symboles correspondant est 
\[S(\R^{2},m)=
\{a\in C^{\infty}(\R^{2}),\,
|\partial^{\alpha}_{\rho}a(\rho)|\le C_{\alpha}m(\rho),\,
\rho\in\R^{2},\alpha\in \N^{2}\}.\]

Nous désignons par  $P=p^{w}$ son $h$-quantifié de Weyl (voir par exemple
\cite{DS}), 
que nous considérons dans $L^{2}(\R).$
Nous faisons une hypothèse d'ellipticité à l'infini sur $P:$ 
\begin{align}\label{int.9}
&|p(x,\xi)|\ge m(x,\xi)/C,\quad|(x,\xi)|\ge C, \nonumber\\
&m(x,\xi)\to \infty,\, (x,\xi)\to \infty.
\end{align}
Pour garantir  que le spectre de $P$ n'est pas le plan complexe, 
nous avons besoin de supposer que 
\begin{equation}\label{int.10}
p(T^{\ast}\R^{2})\ne \C.
\end{equation}
Dans ces conditions, le spectre de $P$ est discret pour $h$ assez petit, voir par 
exemple \cite{Hager2} ou \cite{HS}. 

Introduisons l'ensemble
\begin{equation}
\Sigma(p)=p( T^{\ast}\R),
\end{equation}
ce dernier est automatiquement fermé à cause de (\ref{int.9}).
Nous allons faire une hypothèse sur un point $z_{0}\in \partial\Sigma(p)$:
\begin{equation}\label{int.13}
 p^{-1}(z_{0})=\{\rho_{0}\},  
 \quad\{p,\{p,\bar{p}\}\}(\rho_{0})\ne0.
 \end{equation}
Ceci entraine que $dp(\rho)\ne 0,$
Nous disons aussi $p$ est de type fini d'ordre $2$ au point  $z_{0}.$

Le champ Hamiltonien  $H_{\frac{1}{2i}\{p,\bar{p}\}}$ est un champ de vecteurs réel
et tangent à l'ensemble $\{\rho|\, \frac{1}{2i}\{p,\bar{p}\}(\rho)=0\}.$
Nous avons 
\[]-T_{1},T_{0}[\ni s\mapsto \rho(s):=\exp(sH_{\frac{1}{2i}\{p,\bar{p}\}} )
(\rho_{0})\subset\{
\frac{1}{2i}\{p,\bar{p}\}=0\}\]
est une courbe orientée. 
Sa courbe image par $p$; $\gamma=p\circ\rho$ est une courbe $C^{\infty}$
avec $\dot{\gamma}\ne0$ et dont l'image coïncide avec $\partial\Sigma$ 
près de $z_{0}$. De plus $\Sigma$ se situe à gauche quand on regarde 
dans la direction de $\dot{\gamma},$ voir \cite{Hager3}. 
Près de $z_{0},$ pour tout complexe $z_{\ast}\in\partial\Sigma(p),$ nous avons
$p^{-1}(z_{\ast})=\{\rho(s_{\ast})\}.$ 
$\gamma=p\circ\rho$ est aussi une courbe orientée. 
Le vecteur unitaire tangent à $\partial \Sigma$ dans un voisinage 
de $z_{0}$ est alors donné par 
$u(s):=\dot{\gamma}(s)/|\dot{\gamma}(s)|,$ lequel vérifie 
\[u(s)=-\frac{\{p,\frac{1}{2i}\{p,\bar{p}\}\}(\rho(s))}
{|\{p,\frac{1}{2i}\{p,\bar{p}\}\}(\rho(s))|}.\]
Près de $z_{0}$, $\Sigma(p)$ est un ensemble à bord $C^{\infty}$.
Il existe un voisinage $W$ de $z_{0}$  tel que tout 
$z\in \Sigma(p)\cap W$ peut s'écrire
sous la forme
\begin{equation}
z=\gamma(s)+i\alpha u(s),\quad \alpha\ge 0,
\quad z_{0}=\gamma(0).
\end{equation}
\begin{prop}
Soit $P$ un opérateur $h$-pseudodifférentiel de symbole indépen\-dant de $h.$
Nous supposons que 
$(p,z_{0})$ vérifie les conditions (\ref{int.9}), (\ref{int.10}),
(\ref{int.13}).  Il existe un voisinage $W$ de $z_{0}$ tel que pour 
 tout point $z$ de $(\Sigma(p)\cap W)\setminus\partial\Sigma(p)$ nous avons 
\begin{equation}
p^{-1}(z)=\{\rho_{+}(z),\rho_{-}(z)\},
\quad \rho_{\pm}=(x_{\pm},\xi_{\pm}),
\end{equation}
avec
\[ 0<\pm\frac{1}{2i}\{p,\bar{p}\}(\rho_{\pm})\asymp \alpha^{1/2}\quad 
|\rho_{+}-\rho_{-}|\asymp\alpha^{1/2}.\]
(le point $\rho_{0}$ se scinde en deux points $\rho_{+}$ et 
$\rho_{-}$ lorsque l'on s'éloigne du bord de $\Sigma(p))$.)
\end{prop}
Ici on introduit $\ell_{0}(z).$
Dans le cas où $p'_{\xi}(\rho_{0})\ne0$ nous avons la factorisation 
\begin{equation}
p(\rho)-z=q(x,\xi,\alpha,s)(\xi+g(x,\alpha,s))
\end{equation}
où si $\alpha<0$ alors $\im g(x)>0$, $x>x_{-}$ ou $x<x_{+}$
et $\im g(x)<0,$ $ x_{-}<x<x_{+}.$ Ici $\rho_{\pm}(x_{\pm},\xi_{\pm}),$
$\xi_{\pm}=-g(x_{\pm},\alpha,s).$
Sous ces conditions, on définit $\ell_{0}(z)$ comme
\begin{equation}
\ell_{0}(z):=-\int_{x_{+}(z)}^{x_{-}(z)}g(x)dx.
\end{equation}
Si $p'_{\xi}(\rho_{0})=0$ alors $p'_{x}(\rho)\ne0.$ On applique 
la discussion précédente à $\widehat{p}(x,\xi)=p(-\xi,x),$ fonction
obtenue de $p$ par composition avec une transformation canonique.
On pose alors 
\begin{equation}
\widehat{\ell}_{0}(z) =-\int_{x_{+}}^{x_{-}}\widehat{g}(x,\xi)dx
\end{equation}
où $\widehat{p}(\rho)-z=\widehat{q}(x,\xi,\alpha,s)(\xi+\widehat{g}(x,\alpha,s))$.
On peut prouver que
\begin{prop} 
Si $p'_{\xi}(\rho_{0})\ne0$ et $p'_{x}(\rho_{0})\ne0,$ alors
nous avons
 $\im \ell_{0}(z)=\im \widehat{\ell}_{0}(z)+\MO(\alpha^{\infty}).$
 \end{prop}

\begin{theo}\label{int3}
Il existe une constante $T_{\ast}$ ($<T_{0},T_{1}$) telle que pour toutes constantes 
$C_{0},C_{1}>0$ il existe une constante $C_{2}>0$ telle que la résolvante 
$(P-z)^{-1}$ est bien définie pour 
\begin{equation}
 |s|<T_{\ast},\;
\frac{ h^{2/3}}{C_{0}}\le \alpha\le C_{1}(h\ln \frac{1}{h})^{2/3},\;
 h<\frac{1}{C_{2}},
 \end{equation}
 et satisfait l'estimation
\begin{align*}
\|(P-z)^{-1}\|
\sim&  \frac{\sqrt{\pi}\exp(\frac{1}{h}\im \ell_{0}(z))}{h^{1/2}
(\frac{1}{2i}\{p,\bar{p}\}(\rho_{+}))^{\frac{1}{4}}
(\frac{1}{2i}\{\bar{p},p\}(\rho_{-}))^{\frac{1}{4}}}
\times(1+\mathcal{O}( \tilde{h}))
+\mathcal{O}(\frac{1}{\sqrt{h}\,\alpha^{1/4}}),
\end{align*}
où $\tilde{h}=h/\alpha^{3/2},$ et $\ell_{0}$ (qui est une intégrale d'action)
vérifie 
\begin{equation}
\ell_{0}(z):= \int_{\gamma\subset p^{-1}(z)}\xi dx,
(\gamma \mbox{ relie } \rho_{-}\mbox{ à }\rho_{+}).
\end{equation}
\end{theo}

\textbf{Remerciements:} 
Mes remerciements vont d'abord à Johannes Sj\"ostrand 
pour m'avoir proposé ce sujet, et les aides et les discussions qu'il m'a apportées
ainsi que pour son soutien.
Je tiens aussi à remercier Bernard Helffer pour les échanges sur l'opérateur d'Airy.   
Ce travail a été soutenu par le Grant à Vienne sponsorisé par K. Groechenig 
et  H.G. Feichtinger et par les grants sponsorisés  l'un par V. Ivrii et l'autre par M. Sigal.
Merci enfin à
K. Groechenig  et  H.G. Feichtinger
 à Vienne, et V. Ivrii  et  M. Sigal.
à Toronto pour leur accueil.
\section{Modèle sur $\R$ à paramètre}
Nous rappelons  $p_{\alpha}(x,\xi):=\xi+g(x,\alpha)$ 
désigne le symbole principal semiclassique de 
$P_{\alpha}=hD_{x}+\widetilde{g}(x,\alpha;h).$
Notre opérateur $hD_{x}+\widetilde{g}(x,\alpha;h)$ 
admet comme domaine naturel l'espace semiclassique
suivant
\[H_{sc}(\R):=\{u\in L^{2}|\;\|u\|_{H_{sc}}:=
\|u\|+
\|x^{2}u\|+\|hD_{x}u\|<\infty\}.\]
%
%
%
%
L'inverse, définie pour tout $\alpha$
\begin{align}
(hD_{x}+\widetilde{g}(x,\alpha;h))^{-1}: L^{2}(\R)\to L^{2}(\R),
\end{align}
est donnée par la formule suivante
\[u(x)=\int_{+\infty}^{x}
\frac{i}{h}e^{-\frac{i}{h}\int_{y}^{x}\widetilde{g}(t,\alpha;h) dt}v(y) dy.
\]
En particulier, quand $\alpha$ est nul, nous avons 
(voir \cite{DSZ})
 \begin{equation}
 \|(hD_{x}+\widetilde{g}(x,0;h))^{-1}\| \le Ch^{-2/3},
 \end{equation}
%
%

\subsection{Ensemble d'énergie pour de petites valeurs du paramètre}
Nous cherchons à décrire l'ensemble 
$p_{\alpha}^{-1}(0)$ pour $\alpha\ll 1$ assez petit.
D'abord remarquons que 
\[ p^{-1}_{0}(0)=\{(0,0)\}.\]

Résoudre $\xi+g(x,\alpha)=0$ revient à trouver les 
solutions de l'équation $\im g(x,\alpha)=0.$

Puisque $\im g'_{\alpha}(0)=-1,$ nous avons par la factorisation 
$\im g(x,\alpha)=q(x,\alpha)(\alpha-f(x)),$ où $q\ne 0$ dans un voisinage de 
zéro, avec $f(0)=f'(0)=0,$
et $f''(0)=-\im g''_{xx}(0)\ne 0.$
Par le théorème de l'inversion locale, il existe
des voisinages $U_{+}:=]0,b_{+}[ $ et $ U_{-}:=] b_{-},0[ $ tels que  
$f$ soit un $C^{1}$ diffémorphisme  de $U_{+}$ sur 
$]0,f(b_{+})[$ et de $U_{-}$ sur $]0,f(b_{-})[.$ Pour avoir les solutions 
de l'équation $\alpha-f(x)=0,$ nous inversons simplement une 
fonction croissante en décroissante. Nous trouvons alors 
\begin{align*}
x\in U_{\pm},\quad 
\frac{2}{\im g''_{xx}(0,0)}\alpha
&=x^{2}+\MO(|x|^{3}),\\
&=(x+\MO(x^{2}))^{2}.
\end{align*}
Soit 
\begin{align}
x\in U_{\pm},\quad
x+\MO(x^{2})&=\pm \alpha^{\frac{1}{2}}\left(\frac{2}{\im g''_{xx}(0,0)}\right)^{\frac{1}{2}},\nonumber\\
x(\alpha)&=\pm\alpha^{\frac{1}{2}}\left(\frac{2}{\im g''_{xx}(0,0)}\right)^{\frac{1}{2}}+\MO(\alpha).
\label{g.0}
\end{align}
Nous déduisons alors le résultat suivant:
\begin{lemme} \label{g1}
Soit $p_{\alpha}=\xi+g(x,\alpha)$ le symbole principal de $P_{\alpha}.$
 (\ref{int.4}) et (\ref{int.5}) sont vérifiés.
Nous avons, pour $\alpha>0,$ assez petit, 
\[p^{-1}_{\alpha}(0)=\{\rho_{+}(\alpha)\}\cup\{\rho_{-}(\alpha)\},
\quad\rho_{\pm}(\alpha)=(x_{\pm},\xi_{\pm}),\]
avec
\begin{equation}\label{g.1}
x_{\pm}(\alpha)=\mp \alpha^{\frac{1}{2}}\left(
\frac{-2}{\{p_{0},\frac{1}{2i}\{p_{0},\overline{p_{0}}\}\}(0,0)}
\right)^{1/2}
+\MO(\alpha),
\end{equation}
tel que
\begin{equation}
\pm\frac{1}{2i}\{p_{\alpha},\overline{p_{\alpha}}\}(\rho_{\pm}(\alpha))>0.
\end{equation}
\end{lemme}
L'ensemble $p_{0}^{-1}(0)$ est constitué d'un 
seul point $(0,0)$ qui va se ``scinder'' en 
deux $\rho_{+}(\alpha)$ et $\rho_{-}(\alpha),$ lorsque 
$\alpha$ ne s'annule plus.\\
\noindent \textbf{Preuve.} Il faut remarquer que pour le premier crochet 
et le second crochet de Poisson, nous avons
\begin{align*}
\frac{1}{2i}\{p_{\alpha},\overline{p_{\alpha}}\}(\rho_{\pm})
&=-\im g'_{x}(x_{\pm}(\alpha),\alpha),\\
\{p_{0},\frac{1}{2i}\{p_{0},\overline{p_{0}}\}\}(0,0)&=-\im g_{xx}''(0,0).
\end{align*}
Dès lors,
\begin{align}
\pm\frac{1}{2i}\{p_{\alpha},\overline{p_{\alpha}}\}(\rho_{\pm})&=\mp\im g'_{x}(x_{\pm}(\alpha),\alpha)
\nonumber\\
&=\mp \left(
x_{\pm}(\alpha)\im g''_{xx}(0,0)+\alpha \im g''_{x,\alpha}(0,0)
+\MO(x_{\pm}^{2}+\alpha^{2})
\right)
\label{g.2}
\\
&>0. \nonumber
\end{align}
Le premier terme de (\ref{g.2}) domine car il se comporte comme $\sqrt{\alpha}$.
\hfill $\square$ \medskip

\subsection{Problème de Grushin.} 
Nous nous appuyons sur les sections 
``solutions locales'' et ``problème de Grushin'' de Hager
\cite{Hager1} et 
``Enoncé et résolution asymptotique du problème de Grushin'' 
de Hager \cite{Hager2}. Pour une lecture ``plus aisée'',
nous reprenons en détail les grandes lignes, et les calculs 
dont nous aurons besoin.

On suppose que $\alpha\ge\frac{1}{\MO(1)}$ 
appartient à un intervalle ouvert $\Omega,$ relativement compact,
séparé de l'origine par une 
constante indépendante de $h.$
Nous allons d'abord considérer cette situation, puis
nous laisserons $\alpha$ devenir petit et nous 
ferons des dilations pour se ramener à ce premier cas.

Avec les hypothèses choisies pour $g,$ nous avons 
\begin{equation}
\forall \alpha\in \Omega,\;p^{-1}_{\alpha}(0)=\{\rho_{+}(\alpha),
\rho_{-}(\alpha)\} \; \mbox{ avec }
\pm\frac{1}{2i}\{p_{\alpha},\overline{p_{\alpha}}\}(\rho_{\pm})>0. 
\end{equation}
%
%

On peut trouver des intervalles  $J_{+}$ et $J_{-}$ disjoints
tels que \[\overline{x_{\pm}(\Omega)}\Subset J_{\pm}.\]
Considérons $\chi_{\pm}\in C^{\infty}(I_{\pm})$ avec 
$I_{+}:=]-\infty, \inf J_{-}[$ et 
$I_{-}:=]\sup J_{+},+\infty[.$ Il est clair que 
$I_{+}\cap I_{-}=]\sup J_{+},\inf J_{-}[. $
On demande que $\chi_{\pm}=1$ sur $\overline{J_{\pm}}$ tels que 
$\mbox{supp}(\chi_{+})\cap \mbox{supp}(\chi_{-})=\emptyset.$
Pour le comportement de $\chi_{\pm}$ à l'infini, nous imposons que
$\chi_{+}=1$ sur $]- \infty,\inf J_{+}[$ et,
$\chi_{-}=1$ sur $]\sup J_{-},+\infty[.$
Les $\chi_{\pm}$ s'annulent près de $\partial I_{\pm}.$ 
 
Considérons sur $I_{+},$ une solution de l'équation 
$P_{\alpha}e_{+}=0,$ 
\begin{equation}\label{g.5}
e_{+}:=c_{+}(\alpha;h)\exp{(-\frac{i}{h}\int_{x_{+}}^{x}\widetilde{g}(y,\alpha;h)dy)}.
\end{equation}
Si nous choisissons le formalisme ``BKW'', $e_{+}$ s'écrit
\begin{align*}
e_{+}&=\left(c_{+}(\alpha;h)e^{-ih\int_{x_{+}}^{x}g_{\ast}(y,\alpha;h)dy}
\right) 
e^{-\frac{i}{h}\int_{x_{+}}^{x}g(y,\alpha)dy},\\
&=:a_{+}(x,\alpha;h)e^{\frac{i}{h}\varphi_{+}(x,\alpha)} ,
\quad \varphi_{+}(x,\alpha):=-\int_{x_{+}}^{x}g(y,\alpha)dy,
\end{align*}
avec une phase $\varphi_{+}\in C^{\infty}$, indépendante de $h$, 
qui vérifie l'équation eikonale 
$\partial_{x}\varphi_{+}(x,\alpha)+g(x,\alpha)=0,$
et une amplitude $a_{+}$ admettant un développement asymptotique 
en puissances de $h$,
\[a_{+}(x,\alpha;h)\sim \sum_{k\ge 0}a_{+,k}(x,\alpha)h^{k}
\mbox{ dans } C^{\infty}_{b}(\R),\; \forall \alpha\in\Omega.\]
Nous avons $\varphi_{+}(x_{+},\alpha)=0,$ $\varphi'_{+}(x_{+},\alpha)=\xi_{+}\in \R,$
et 
\[\im \varphi''_{+}(x_{+}(\alpha),\alpha)=
-\im g'_{x}(x_{+},\alpha)=
\frac{1}{2i}\{p_{\alpha},\overline{p_{\alpha}}\}>0.\]
La partie imaginaire de la phase est positive, cela restera
vrai globalement sur $I_{+}:$ $\im \varphi_{+}\ge 0.$

Grâce au choix de $g$ vers $-\infty,$
pour $C>0$ assez grand nous avons 
$\|e_{+}\|_{L^{2}(]-\infty,-C])}=\MO(e^{-\frac{1}{\widetilde{C}h}}).$
Il est alors possible, grâce au lemme de la phase stationnaire
(que nous rappellerons après), de choisir $c_{+}(\alpha;h)$ de la forme 
\[c_{+}(\alpha;h)\sim h^{-1/4}(c_{+}^{0}(\alpha)+ hc^{1}_{+}(\alpha)+\ldots)>0,\]
avec
\begin{align*}
c_{+}^{0}(\alpha)&= 
\left(\frac{ -2 \im g'_{x}(x_{+}(\alpha),\alpha)^{1/4}}{2\pi}\right)^{\frac{1}{4}}
\\
&= \frac{(- \im g'_{x}(x_{+},\alpha))^{1/4}}{\pi^{1/4}},
\end{align*}
tel que  $e_{+}$ soit normalisé dans $L^{2}$  sur $I_{+}.$
Ce qui implique, dans le formalisme ``BKW'', que 
\begin{align*}
a_{+}(x,\alpha;h)&=c_{+}(\alpha;h)e^{-ih\int_{x_{+}}^{x}g_{\ast}(y,\alpha;h)dy}\\
&=c_{+}(\alpha;h) (1-ih\int_{x_{+}}^{x}g^{0}_{\ast}(y,\alpha)dy+\MO(h^{2})),\\
&= \frac{| \im g'_{x}(x_{+},\alpha)|^{1/4}}{(\pi h)^{1/4}}(1+\MO(h))
\mbox{ dans } C^{\infty}_{b}(\R).
\end{align*}
\begin{lemme}\label{g11} On désigne par $e_{+}\in H_{sc}(I_{+})$ la solution normalisée 
dans $L^{2}(I_{+})$ de $(hD_{x}+\widetilde{g}(x,\alpha;h))e_{+}=0$ sur $I_{+}.$
Nous avons la représentation asymptotique suivante pour 
$h\ll 1,$ assez petit 
 \[
 e_{+}\sim \frac{ |\im g'_{x}(x_{+},\alpha)|^{1/4}}{(\pi h)^{1/4}}(1+\MO(h))
e^{-\frac{i}{h}\int_{x_{+}}^{x}\widetilde{g}(y,\alpha;h)dy}.\]
\end{lemme}

De manière analogue nous avons :
\begin{lemme}\label{g12} 
On désigne par $e_{-}\in H_{sc}(I_{-})$ 
la solution normalisée dans 
$L^{2}(I_{-})$ de $(hD_{x}+\widetilde{g}(x,\alpha;h))^{\ast}e_{-}=0$ sur $I_{-}.$
$e_{-}$ admet la représentation asymptotique suivante
 pour $h\ll 1,$ assez petit 
 \[e_{-}\sim \frac{(\im g'_{x}(x_{-},\alpha))^{1/4}}{(\pi h)^{1/4}}(1+\MO(h))
e^{-\frac{i}{h}\int_{x_{-}}^{x}\overline{\widetilde{g}(y,\alpha;h)}dy}.\]
\end{lemme}

Nous rappelons ici le lemme de la phase stationnaire:
\begin{prop}[Phase stationnaire]
Soit $a\in C^{\infty}_{0}(\R).$ Supposons que $0\in K=\mathrm{supp}(a)$
et 
\[\varphi(0)=\varphi'(0)=0,\quad\varphi''(0)\ne 0,\quad
\im \varphi\ge 0.\]
Supposons en plus que $\varphi'(x)\ne 0$ sur $K-\{0\}.$
Posons $g(x)=\varphi(x)-\varphi''(0)x^{2}/2.$
Alors, nous avons le développement asymptotique explicite, lorsque 
$h \to 0,$ 
\begin{align}
&\int_{\R}e^{i\varphi(x)/h}a(x)dx\sim\nonumber\\
\left(\frac{2\pi ih}{\varphi''(0)}\right)^{\frac{1}{2}}
\sum_{k\ge 0}\sum_{\ell \ge 0}&\left(
\frac{h}{2i\varphi''(0)}\right)^{k}\frac{1}{\ell !}\frac{1}{k!}\frac{d^{2k}}{dx^{2k}}
((i/h)^{\ell}g^{\ell}a)(0). 
\end{align}
En particulier, nous voyons que le premier terme est 
$(2\pi ih)^{1/2}\varphi''(0)^{-1/2}a(0).$
\end{prop}
Pour une preuve on pourra consulter \cite{ZW} p.36.
\begin{prop}\label{g13}
Pour $v\in L^{2}(I_{+}),$ $v_{+}\in \C,$ le problème de Cauchy 
$(hD_{x}+\widetilde{g}(x,\alpha;h))u=v, u(x_{+})=0,$
admet la solution unique $u=\widetilde{F}v$ avec 
\[ \|\widetilde{F}\|_{L^{2}(I_{+})\to H_{sc}(I_{+})}\le \frac{C}{\sqrt{h}}.\]
\end{prop}
Par ailleurs $\widetilde{F}$ admet un noyau intégral $k$ vérifiant l'égalité
\[|k(x,y)|\le\MO(1/h)\exp(-|x-y|/(C\sqrt{h})).\]
\textbf{Preuve.}
$\widetilde{F}$ admet le noyau intégral
\begin{equation}
k(x,y)=
\frac{i}{h} e^{-\frac{i}{h}\int^{x}_{y}\widetilde{g}(\tilde{x},\alpha;h) \,d\tilde{x}}
1_{\{x_{+}\le y\le x\}}
-
\frac{i}{h} e^{-\frac{i}{h}\int^{x}_{y}\widetilde{g}(\tilde{x},\alpha;h) \,d\tilde{x}}
1_{\{x\le y\le x_{+}\}}.
\end{equation}
Nous montrerons après qu'un tel noyau vérifie l'inégalité suivante
\begin{equation}\label{g.5bis}
|k(x,y)|\le\MO(1/h)\exp(-|x-y|/(C\sqrt{h}))
\end{equation}
pour une constante $C$ donnée.

La norme $L^{2}$ de $\widetilde{F}$  est majorée (lemme de Schur) par
\begin{equation}\label{g.6}
\left( \sup_{x}\int |k(x,y)|dy\right)^{1/2}
\left( \sup_{y}\int |k(x,y)|dx\right)^{1/2}.
\end{equation}
Ce qui nous donne
\begin{align}
\sup_{x\in I_{+}}
\int |k(x,y)| dy\le \frac{C}{\sqrt{h}}\mbox{ et, }
\sup_{y\in I_{+}}\int|k(x,y)|dx\le \frac{C}{\sqrt{y}}.
\end{align}
Donc
\[\|\widetilde{F}\|_{L^{2}(I_{+})\to L^{2}(I_{+})}\le \frac{C}{\sqrt{h}}.\]
De la même fa\c con on peut montrer que
\[\|x^{2}\widetilde{F}\|_{L^{2}(I_{+})\to L^{2}(I_{+})}\le \frac{C}{\sqrt{h}}.\]
Ce qui implique, en utilisant 
$(hD_{x}+\widetilde{g}(x,\alpha;h))\widetilde{F}w=w,$
que nous avons
\[\|\widetilde{F}\|_{L^{2}(I_{+})\to H_{sc}(I_{+})}\le \frac{C}{\sqrt{h}}.\]
Il nous reste à prouver l'estimation (\ref{g.5bis}).
Parce que $x_{+}$ est un zéro non-dégénéré de $\im g,$
et au comportement de $\im g$ vers $-\infty,$ nous avons
pour $h$ assez petit, et $x\in I_{+},$ (rappelons que 
$g_{\ast}(x,\alpha;h)$ est uniformément borné)
\begin{align}\label{g.7}
\im \widetilde{g}(x,\alpha;h)&\le -\frac{1}{C_{0}}(x-x_{+})
+C_{0}h^{2} \mbox{ si } x\ge x_{+},\\
\im \widetilde{g}(x,\alpha;h)&\ge-\frac{1}{C_{1}}(x-x_{+}) 
-C_{1}h^{2}\mbox{ si } x\le x_{+}.
\label{g.8}
\end{align}
Ce qui implique, si nous considérons le cas $x\ge x_{+},$
\begin{align}
|k(x,y)| &\le \frac{1}{h} e^{\frac{1}{h}\int^{x}_{y}\im 
\widetilde{g}(\tilde{x},\alpha;h) \,d\tilde{x}}
1_{\{x_{+}\le y\le x\}}
\nonumber\\
&\le
\frac{1}{h} e^{-\frac{1}{C_{0}h}\int^{x}_{y}
(\widetilde{x}-x_{+})\,d\tilde{x}}
e^{C_{0}h(x-y)}
1_{\{x_{+}\le y\le x\}}
\nonumber\\
&\le \frac{1}{h} e^{-\frac{1}{2C_{0}h}((x-x_{+})^{2}-(y-x_{+})^{2})}
e^{C_{0}h(x-y)}
1_{\{x_{+}\le y\le x\}},\label{g.8bis}
\end{align}
Pour $(x-x_{+})\ge \sqrt{h},$ \ref{g.8bis} s'écrit 
\[\frac{1}{h} e^{-\frac{1}{2C_{0}h}((x-x_{+})-2(C_{0}h)^{2})(x-y)}1_{\{x_{+}\le y\le x\}}=
\frac{1}{h} e^{-|x-y|/(C\sqrt{h})},
\]
et pour $(x-x_{+})\le \sqrt{h},$
\[\frac{1}{h} e^{-\frac{1}{2C_{0}h}((x-x_{+})-2(C_{0}h)^{2})(x-y)}1_{\{x_{+}\le y\le x\}}=
\MO(1/h) e^{-|x-y|/(C\sqrt{h})}.\]
Les autres inégalités nécessaires pour déduire (\ref{g.5bis}) se calculent de la même manière.
\hfill $\square$ \medskip

\begin{prop}
Pour $P_{\alpha}=hD_{x}+\widetilde{g}(x,\alpha;h),$ $v\in L^{2}(I_{+}),$ $v\in \C$ le problème 
\begin{equation}
\left\{
\begin{array}{cc}
P_{\alpha}u=v\\
R_{+}u=v_{+}
\end{array}
\right.
\end{equation}
avec 
\[R_{+}u:=\langle u,\chi_{+}e_{+}\rangle =\int_{I_{+}}u(x)\chi_{+}(x)\overline{e_{+}(x)}dx,\]
admet une solution unique 
\[ u=Fv+F_{+}v_{+}\in H_{sc},\]
où 
\[F_{+}v_{+}:=\frac{1}{\langle e_{+},\chi_{+}e_{+}\rangle}v_{+}e_{+}
=:\frac{1}{D_{+}}v_{+}e_{+},\]
avec 
\[D_{+}=1+\MO(e^{-\frac{1}{Ch}}),\]
et nous avons 
\begin{align}
\|F\|_{L^{2}\to H_{sc}}&=\MO(h^{-1/2}),\\
\|F\|_{\C\to H_{sc}}&=\MO(1).
\end{align}
\end{prop}
\textbf{Preuve.} Le problème 
\[\left\{
\begin{array}{cc}
P_{\alpha}u=0\\
R_{+}u=v_{+}
\end{array}
\right.\]
admet la solution unique $u=F_{+}v_{+},$
tandis que le problème 
\[\left\{
\begin{array}{cc}
P_{\alpha}u=v\\
R_{+}u=0
\end{array}
\right.\]
admet la solution unique $u=Fv:=(1-F_{+}R_{+})\widetilde{F}v$.
\hfill $\square$ \medskip

Nous désignons par  $L^{2}_{comp}(I_{-})$ et $H_{sc,comp}(I_{-})$
respectivement l'ensemble des fonctions $f$ appartenant à $L^{2}$ 
et à $H_{sc}(I_{-})$ respectivement telles que $f$ est nulle 
dans un voisinage de $\inf I_{-}.$

\begin{prop}
Pour $P_{\alpha}=hD_{x}+\widetilde{g}(x,\alpha;h),$ $v\in L^{2}_{comp}(I_{-})$ le problème
\[P_{\alpha}u+R_{-}u_{-}=v,\]
avec 
\[R_{-}u_{-}:=u_{-}\chi_{-}e_{-},\;
u_{-}\in \C,\]
admet une solution unique dans $H_{sc,comp}\times\C,$ donnée par
\begin{align}
u&=Gv,\\
u_{-}&=G_{-}v:= \frac{1}{\langle \chi_{-}e_{-},e_{-}\rangle }
\langle v,e_{-}\rangle
=:\frac{1}{D_{-}}\langle v,e_{-}\rangle,
\end{align}
où
\[ D_{-}=1+\MO(e^{-\frac{1}{Ch}}),\]
avec, en tant qu'opérateur $L^{2}_{comp}(I_{-})\to H_{sc,comp}(I_{-}):$
\[\|G\|_{L^{2}\to H_{sc}}\le \frac{C}{\sqrt{h}}. \]
\end{prop}
\textbf{Preuve.}
Pour $x\le x_{-},$ nous avons l'unique solution de 
$P_{\alpha}u=\tilde{v},$ $\tilde{v}\in L^{2}_{comp}(I_{-}),$
qui est nulle près de $\inf I_{-}$
\[u_{1}(x)=\frac{i}{h}
\int_{-\infty}^{x}e^{-\frac{i}{h}\int^{x}_{y}\widetilde{g}(\tilde{x},\alpha;h)d\tilde{x}}
\tilde{v}(y)dy:=\widetilde{G}_{1}\tilde{v}(x),\]
alors que pour $x\ge x_{-},$ nous avons 
\[u_{2}(x)=\frac{i}{h}
\int_{+\infty}^{x}e^{-\frac{i}{h}\int^{x}_{y}\widetilde{g}(\tilde{x},\alpha;h)d\tilde{x}}
\tilde{v}(y)dy:=\widetilde{G}_{2}\tilde{v}(x),\]
les deux exposants étant décroissants (à partie réelle strictement négative
loin de $x_{-}$) dans le domaine d'intégration.  
Afin d'obtenir une solution continue, il faut imposer 
\begin{align*}
0=u_{1}(x_{-})-u_{2}(x_{-})&=\frac{i}{h}\int^{+\infty}_{-\infty}
e^{-\frac{i}{h}\int^{x_{-}}_{y}\widetilde{g}(\tilde{x},\alpha;h)d\tilde{x}}\tilde{v}(y)dy,\\
&=\frac{i}{h\overline{c_{-}}}\langle \tilde{v},e_{-}\rangle.
\end{align*}
Donc, avec 
\[u_{-}=G_{-}v,\]
nous pouvons prendre $\tilde{v}=v-R_{-}u_{-},$ car par construction
\[\langle v-R_{-}u_{-},e_{-}\rangle=0.\]
Nous avons alors la solution 
\begin{align*}
u&=Gv=\widetilde{G}(I-R_{-}G_{-})v,\\
u_{-}&=G_{-}v,
\end{align*}
où $\widetilde{G}$ est donné par $\widetilde{G}_{1,2}$ dans les zones correspondantes.
En observant que $\widetilde{G}$ a le noyau intégral
\[
k(x,y)=
\frac{i}{h} e^{-\frac{i}{h}\int^{x}_{y}\widetilde{g}(\tilde{x},\alpha;h) \,d\tilde{x}}
1_{\{y\le x\le x_{-}\}}
-
\frac{i}{h} e^{-\frac{i}{h}\int^{x}_{y}\widetilde{g}(\tilde{x},\alpha;h) \,d\tilde{y}}
1_{\{x_{-}\le x\le y\}},
\]
- qui est semblable au noyau intégral de $\widetilde{F}$- et que 
$\im g'_{x}(x_{+})\sim -\im g'_{x}(x_{-}),$ nous pouvons utiliser les estimations 
du paragraphe précédent pour trouver que   
\[\|\widetilde{G}\|_{L^{2}(I_{-})\to L^{2}(I_{-})}\le \frac{C}{\sqrt{h}},\]
ainsi que
\[\|G\|_{L^{2}(I_{-})\to H_{sc}(I_{-})}\le \frac{C}{\sqrt{h}}.\]
\hfill $\square$ \medskip

\paragraph{Problème de Grushin global.}
Nous choisissons une partition de l'unité de $\R,$
$\psi_{\pm}\in C^{\infty}(I_{\pm}),$  $ \psi_{+}+\psi_{-}=1,$
telle que
$\chi_{\pm} \prec \psi_{\pm},$ où
$\psi\prec \phi$ signifie
\begin{align}
\mbox{supp}\,(\psi)\cap\mbox{supp}(1-\phi)=\emptyset.
\end{align}
Comme les $\chi_{\pm},$ les $\psi_{\pm}$ s'annulent près du bord de l'intervalle $I_{\pm}$.
On suit Hager \cite{Hager1} ``Inverse global'' :
on commence par résoudre le problème sur $I_{+}$
\[
\left\{
\begin{array}{cc}
(hD_{x}+\widetilde{g}(x,\alpha;h))u=\psi_{+}v\\
R_{+}u=v_{+}.
\end{array}
\right.
\]
Sur $I_{+},$ on introduit
\[u_{1}:=(1-\chi_{-})F\psi_{+}v+(1-\chi_{-})F_{+}v_{+}.\]
Donc, en utilisant $\chi_{+}\prec(1-\chi_{-}),$ et 
$\psi_{+}\prec(1-\chi_{-})$ (car $\chi_{-}\prec\psi_{-})$
\begin{align*}
&R_{+}u_{1}=\langle u_{1},\chi_{+}e_{+}\rangle =v_{+},\\
&P_{\alpha}u_{1}=\psi_{+}v-[P_{\alpha},\chi_{-}]F\psi_{+}v-
[P_{\alpha},\chi_{-}]F_{+}v_{+},
\end{align*}
en rappelant que $P_{\alpha}:=hD_{x}+\widetilde{g}(x,\alpha;h).$

Comme chez Hager \cite{Hager1}, on peut ``corriger l'erreur'' sur $I_{-}$
en y résolvant le problème 
\[(hD_{x}+\widetilde{g}(x,\alpha;h))u_{2}+R_{-}u_{-}=\psi_{-}v
+[P_{\alpha},\chi_{-}]F\psi_{+}v
+[P_{\alpha},\chi_{-}]F_{+}v_{+}.\]
Nous avons la solution
\begin{align*}
 &u_{2}=G(\psi_{-}v
+[P_{\alpha},\chi_{-}]F\psi_{+}v-
[P_{\alpha},\chi_{-}]F_{+}v_{+})\\
&u_{-}=G_{-}(\psi_{-}v
+[P_{\alpha},\chi_{-}]F\psi_{+}v-
[P_{\alpha},\chi_{-}]F_{+}v_{+}).
\end{align*}
Et $R_{+}u_{2}=0$ car $\mbox{supp}(u_{2})
\cap\mbox{supp}(\chi_{+})=\emptyset.$

Finalement nous obtenons :
\begin{theo}\label{g16}
Pour tout  $\alpha\in \Omega,$
\[\mathcal{P}=
\left(\begin{array}{cc}
hD_{x}+\widetilde{g}(x,\alpha;h)& R_{-}\\
R_{+}&0
\end{array}
\right): H_{sc}\times \C\to L^{2}\times \C, \]
est inversible d'inverse 
\[\mathcal{E}=
\left(\begin{array}{cc}
E& E_{+}\\
E_{-}&E_{-+}
\end{array}
\right)\]
où
\begin{align}\label{g.48}
&E=G\left(\psi_{-}+\frac{h}{i}\chi'_{-}F\psi_{+}\right)+ (1-\chi_{-})F\psi_{+} \\
&E_{+}=(1-\chi_{-})F_{+}+G\frac{h}{i}\chi'_{-}F_{+} \nonumber\\
&E_{-}=G_{-}(\psi_{-}+\frac{h}{i}\chi'_{-}F\psi_{+}) \nonumber\\
&E_{-+}=G_{-}\frac{h}{i}\chi'_{-}F_{+}=
-\frac{h}{iD_{-}D_{+}}\langle \chi'_{-} e_{+},e_{-}\rangle. \nonumber
\end{align}
Les normes vérifient
\begin{align}\label{g.49}
&\|E\|_{L^{2}\to H^{1}_{sc}}=\mathcal{O}(\frac{1}{\sqrt{h}}),\,
\|E_{+}\|=\mathcal{O}(1),\,
\|E_{-}\|=\mathcal{O}(1),\\
&\|E_{-+}\|=\mathcal{O}\left(\sqrt{h}\,e^{-\frac{1}{Ch}}\right).
\label{g.50}
\end{align}
Les opérateurs $F,F_{+},G,G_{-} $ ont été définis au-dessus.
De plus $E_{+}$ et $E_{-}$ sont des opérateurs de rang 1, satisfaisant
\begin{align}
E_{+}&=\frac{1-\chi_{-}}{D_{+}} e_{+}+\MO(h^{\infty}),\\
E_{-}&=\langle\bullet ,\frac{\psi_{-}}{D_{-}}e_{-}\rangle+\MO(h^{\infty}).
\end{align}
\end{theo}
On peut déduire une formule exacte pour $E_{-+},$
\begin{equation}
E_{-+}=\frac{h c_{+} \overline{c_{-}}}{iD_{+}D_{-}} 
\exp{(\frac{i}{h}\int_{x_{-}}^{x_{+}}\widetilde{g}(\tilde{x},\alpha;h)d\tilde{x}}).
\end{equation}

Finalement, en utilisant la relation
$(hD_{x}+\widetilde{g})^{-1}=E-E_{+}(E_{-+})^{-1}E_{-}$,
on trouve l'estimation suivante pour la résolvante:
\begin{prop}\label{g18}
On rappelle que $1/\MO(1)\ge\alpha\le\MO(1)$ appartient à un ouvert $\Omega$, relativement compact
et séparé de l'origine par une constante indépendante de h.
Pour tout $\alpha$ dans $\Omega,$ la résolvante vérifie pour $h$ assez petit  
\begin{align}\nonumber
\|(hD_{x}+\widetilde{g}(x,\alpha;h))^{-1}\|\sim &
\frac{\sqrt{\pi}e^{\frac{1}{h}\im \ell_{0}(\alpha)}}{h^{1/2}(\im g'_{x}(x_{+},\alpha))^{\frac{1}{4}}
(\im \bar{g}'_{x}(x_{-},\alpha))^{\frac{1}{4}}}(1+\mathcal{O}(h))\\
&+\mathcal{O}(\frac{1}{\sqrt{h}}),\label{g.52}
\end{align} 
où $\ell_{0}$ vérifie 
\begin{equation}
\ell_{0}:=-\im \int_{x_{+}}^{x_{-}}g(y,\alpha)dy.
\end{equation}
\end{prop}

\textbf{Preuve.} Il suffit de remarquer que $E_{-+}$ s'écrit
\begin{align*}
|E_{-+}(\alpha)| \sim \frac{h^{\frac{1}{2}}} {\sqrt{\pi}}  
\big(\im &g'_{x}(x_{+},\alpha)\,
\im \bar{g}'_{x}(x_{-},\alpha)
\big)^{\frac{1}{4}}\\
&
\exp{(-\frac{1}{h}\im \int^{x_{+}}_{x_{-}}g(y,\alpha)dy)}
\; (1+\mathcal{O}(h)),
\end{align*}
et que 
\begin{equation*}
\|E_{-}E_{+}\|=\|\frac{1-\chi_{-}}{D_{+}}e_{+}\|\,
\|\frac{\psi_{-}}{D_{-}}e_{-}\|
\sim 1+\MO(h^{\infty}).
\end{equation*}
\hfill $\square$ \medskip

\paragraph{L'opérateur d'Airy complexe sur $\R.$} 
Considérons l'opérateur $\mathcal{A}=D^{2}_{x}+ix$ sur la droite réelle,
muni du domaine 
$D(\mathcal{A})=\{u\in H^{2}(\R)|\; xu\in L^{2}(\R)\}.$
Nous pouvons montrer que pour chaque $z\in \C,$ la résolvante satisfait
(voir Helffer \cite{He})
\[\|(\mathcal{A}-z)^{-1}\|=\|(\mathcal{A}-\re z)^{-1}\|.\]
\begin{prop}[J.~Martinet]
Pour $\re z>0$ assez grand, il existe deux constantes positives $C_{0}$
et $C_{1}$ telles que 
\[C_{0}|\re z|^{-\frac{1}{4}}\exp\frac{4}{3}(\re z)^{\frac{3}{2}}\le
\|(\mathcal{A}-z)^{-1}\|\le
C_{1}|\re z|^{-\frac{1}{4}}\exp\frac{4}{3}(\re z)^{\frac{3}{2}}.\]
\end{prop}
nous voyons que 
\[\mathcal{F}(D_{x}^{2}+ix-z)\mathcal{F}^{-1}=i(-D_{x}-ix^{2}+iz).\]
Après conjugaison par l'opérateur unitaire $u(x)\mapsto u(-x),$ 
nous notons $(\mathcal{A}_{0}+iz):=D_{x}-ix^{2}+iz$ ce nouvel opérateur. 
Nous avons 
$\|(\mathcal{A}_{0}-iz)^{-1}\|=\|(\mathcal{A}_{0}-i\re z)^{-1}\|.$
Le changement de variable $x=(\re z)^{\frac{1}{2}}y,$ nous permet 
d'identifier 
$\mathcal{A}_{0}+iz$ avec $|\re z|\mathcal{Q},$ où $\mathcal{Q}$ est donné par 
($y$ a été remplacé par $x$)
\[\frac{1}{(\re z)^{\frac{3}{2}}}D_{x}-ix^{2}+i.\]
Nous nous retrouvons alors dans la situation précédente avec $\im g(x,\alpha)\le0,$
plus explicitement $\mathcal{Q}$ est de la forme $hD_{x}+\widetilde{g}$ avec 
$\widetilde{g}=g.$
En reprenant les notations précédentes, nous avons
\begin{align}
 &h=\frac{1}{(\re z)^{\frac{3}{2}}},\\
&g(x,z):=-ix^{2}+i,\\
&x_{\pm}(z)=\pm1\mbox{ vérifiant } \im g(x_{\pm},z)=0,\\
&-\im g'_{x}(x_{+}(z),z)=2x_{+}(z)=2,\\
&+\im g'_{x}(x_{-}(z),z)=-2x_{-}(z)=2.
\end{align}
De plus
\begin{align}
-\frac{1}{h}\im\int_{x_{+}}^{x_{-}}g(x,z)dx=
-\frac{1}{h}\im\int_{1}^{-1}g(x,z)dx&=
-(-2+\frac{2}{3})(\re z)^{3/2}\nonumber\\
&=\frac{4}{3}(\re z)^{3/2},
\end{align}
et 
\[h^{1/2}\big(\im g'_{x}(x_{+},z)\,
\im \bar{g}'_{x}(x_{-},z)
\big)^{\frac{1}{4}}=
h^{1/2}(2x_{+}(z))^{1/2}=2^{1/2}(\re z)^{-\frac{3}{4}}.\]
On déduit alors le corollaire \ref{int2}.

\subsection{Estimation de résolvante près du bord de l'image du symbole.}
On suppose cette fois ci que $\alpha>0$ tend vers zéro pour le cas modèle à paramètre
$hD_{x}+\widetilde{g}(x,\alpha;h)$ sur la ligne droite.
On va se ramener au cas précédent par ``scaling''.
Rappelons que (lemme \ref{g1})
\begin{equation}
x_{\pm}(\alpha)=\mp\alpha^{1/2}
\left(\frac{-2}{\{p_{0},\frac{1}{2i}\{p_{0},\overline{p_{0}}\}\}(0,0)}\right)^{1/2}+
\MO(\alpha).
\end{equation}
Le changement de variable $x=(\alpha)^{1/2} y,$ nous permet d'identifier 
notre opérateur avec $\alpha Q$ où $Q$ est donné par 
\begin{align}
\frac{h}{\alpha^{1/2}\alpha}D_{y}+\frac{g(\alpha^{1/2}y,\alpha)}{\alpha}
+\frac{h^{2}}{\alpha^{3}}(\alpha^{3}g_{\ast}(\alpha^{1/2}y,\alpha;h)=
\nonumber\\
\tilde{h}D_{y}+\frac{g(\alpha^{1/2}y,\alpha)}{\alpha}
+\tilde{h}^{2}(\alpha^{3}g_{\ast}(\alpha^{1/2}y,\alpha;h))
\mbox{ où } \tilde{h}=h/\alpha^{3/2}.
\end{align}
%
%
On a ainsi remplacé les points $x_{+}(\alpha)$ et $x_{-}(\alpha),$ par des points
$y_{+}(\alpha)=\alpha^{-\frac{1}{2}}x_{+}$ et $y_{-}(\alpha)=\alpha^{-\frac{1}{2}}x_{-},$ 
où $y_{\pm}\asymp \mp 1.$
Les points $y_{\pm}$ vérifient $\im f(y_{\pm},\alpha)=0,$
où $f$ est défini par 
$f(x,\alpha):=\frac{g(\alpha^{1/2}x,\alpha)}{\alpha}.$
Nous avons $\im f''_{xx}(0,0)>0,$
et
$\im f(x,\alpha)\asymp (x^{2}-1),$ 
Nous pouvons alors appliquer 
la proposition \ref{g18} avec 
$\tilde{h}:=h/\alpha^{3/2}$ comme notre 
nouveau petit paramètre semiclassique.
\begin{prop}
Pour tout $\alpha,$ la résolvante vérifie pour $h\ll \alpha^{3/2}$ et $\alpha$ assez petit  
\begin{align}\label{g.59}
\|(hD_{y}+\widetilde{g}(y,\alpha;h)^{-1}\|\sim &
\frac{\sqrt{\pi}
\exp(\frac{1}{\tilde{h}}\im \widetilde{\ell}(\alpha))}
{\alpha\tilde{h}^{1/2}
(\im f'_{y}(y_{+},\alpha))^{1/4}
(\im \overline{ f'_{y}(y_{-},\alpha)})^{1/4}}
(1+\mathcal{O}( \tilde{h}))
\nonumber\\
&+\mathcal{O}(\frac{1}{\alpha\sqrt{\tilde{h}}}),
\end{align} 
où $\tilde{h}=h/\alpha^{3/2},$ et $\widetilde{\ell}$ vérifie 
\begin{equation}
\widetilde{\ell}:=\int_{y_{-}}^{y_{+}}\frac{g(\alpha^{1/2}\tilde{y},\alpha)}{\alpha} d\tilde{y}, \quad 
\im f(y_{\pm},\alpha)=0.
\end{equation}
Rappelons que $f=\frac{g(\alpha^{1/2}y,\alpha)}{\alpha}.$
\end{prop}

Ensuite, nous avons les égalités suivantes, 
\begin{align*}
\frac{1}{\tilde{h}}\widetilde{\ell}&=\frac{\alpha^{3/2}}{h}\int_{y_{-}}^{y_{+}}
\frac{g(\alpha^{1/2}\tilde{y},\alpha)}{\alpha} d\tilde{y}
=\frac{\alpha^{3/2}}{h}\int_{\alpha^{1/2}y_{-}}^{\alpha^{1/2}y_{+}}
\frac{g(\tilde{y},\alpha)}{\alpha\alpha^{1/2}} d\tilde{y}\\
&= \frac{1}{h}\int_{x_{-}}^{x_{+}}g(\tilde{y},\alpha) d\tilde{y}=:
\frac{1}{h}\ell_{0}(\alpha),
\end{align*}
et, 
\begin{align*}
\alpha\tilde{h}^{1/2}
&(-\im f'_{y}(y_{+},\alpha))^{1/4}
\im ( f'_{y}(y_{-},\alpha))^{1/4}\\
&=
\alpha^{1/4}h^{1/2}
\frac{(\im g'_{x}(\alpha^{1/2}y_{+},\alpha))^{\frac{1}{4}}}{\alpha^{1/8}} 
\frac{\im \overline{({g}'_{x}(\alpha^{1/2}y_{-},\alpha)})^{\frac{1}{4}}}{\alpha^{1/8}} \\
&=h^{1/2}
(\im g'_{x}(x_{+},\alpha))^{\frac{1}{4}}
(\im \overline{{g}'_{x}(x_{-},\alpha)})^{\frac{1}{4}}\\
&=
h^{1/2}(\frac{1}{2i}\{p_{\alpha},\overline{p_{\alpha}}\}(\rho_{+}))^{\frac{1}{4}}
(\frac{1}{2i}\{\overline{p_{\alpha}},p_{\alpha}\}(\rho_{-}))^{\frac{1}{4}}.
\end{align*}

Nous pouvons ainsi donner une formule invariante pour la résolvante.
\begin{corol}\label{g22} Soit $g$ satisfaisant (\ref{int.4}) et (\ref{int.5}). 
Pour tout $\alpha,$ la résolvante vérifie pour $h\ll \alpha^{3/2}$ avec $\alpha$ assez petit  
\begin{align}
\|(hD_{x}+\widetilde{g}(x,\alpha;h))^{-1}\|\sim  \frac{\sqrt{\pi}\exp(\frac{1}{h}\im \ell_{0}(\alpha))}{h^{1/2}
(\frac{1}{2i}\{p_{\alpha},\overline{p_{\alpha}}\}(\rho_{+}))^{\frac{1}{4}}
(\frac{1}{2i}\{\overline{p_{\alpha}},p_{\alpha}\}(\rho_{-}))^{\frac{1}{4}}}
(1+\mathcal{O}( \tilde{h}))\nonumber\\
+\mathcal{O}(\frac{1}{\sqrt{h}\,\alpha^{1/4}}),
\end{align} 
où $\tilde{h}=h/\alpha^{3/2},$ et $\ell_{0}$ vérifie 
\begin{equation}
\ell_{0}:=\im \int_{x_{-}}^{x_{+}}g(y,\alpha)dy, \quad 
\im g(x_{\pm},\alpha)=0.
\end{equation}
\end{corol}
\begin{remarque}
Puisque $|\alpha_{\pm}(\alpha)|\asymp \alpha^{\frac{1}{2}},$ 
nous avons les estimations suivantes :
\begin{align*}
&\ell_{0}(\alpha)\asymp \alpha\alpha^{\frac{1}{2}}=\alpha^{3/2},\\
&\frac{1}{2i}\{p_{\alpha},\overline{p_{\alpha}}\}(\rho_{+})=
-\im g'_{x}(x_{+},\alpha)\asymp \alpha^{1/2},\\
&\frac{1}{2i}\{\overline{p_{\alpha}},p_{\alpha}\}(\rho_{-})=
-\im \bar{g}'_{x}(x_{-},\alpha)\asymp \alpha^{1/2}.
\end{align*}
Donc, sous la condition $h\ll \alpha^{3/2}\ll 1$ assez petit, nous avons
\[ \|(hD_{x}+\widetilde{g}(x,\alpha;h))^{-1}\|\le \frac{C}{h^{1/2}\alpha^{1/4}}
e^{\frac{C\alpha^{3/2}}{h}}.\]
Nous retrouvons alors la formule usuelle donnée dans \cite{WB} et \cite{SJ}.
\end{remarque}

Résumons le cas où $\alpha$ tend vers zéro, le théorème \ref{g16}
traite en effet du point $\alpha\ge\frac{1}{\MO(1)}$ :

\begin{theo}\label{g23}
Rappelons que $(hD_{x}-\widetilde{g}(x,\alpha;h))^{-1}=
E-E_{+}E_{-+}^{-1}E_{-},$
avec
\begin{align}\label{g.60}
|E_{-+}^{-1}|\sim
 \frac{\sqrt{\pi}\exp(\frac{1}{h}\im \ell_{0}(\alpha))}{h^{1/2}
(\frac{1}{2i}\{\xi+g,\xi+\bar{g}\}(\rho_{+}))^{\frac{1}{4}}
(\frac{1}{2i}\{\xi+\bar{g},\xi+g\}(\rho_{-}))^{\frac{1}{4}}}
\times(1+\mathcal{O}( \tilde{h}))
\end{align}
$\tilde{h},$
$\ell_{0}$ seront rappelés par la suite, et 
$E=\mathcal{O}(\frac{1}{\sqrt{h}\,\alpha^{1/4}}).$
$E_{+}$ et $E_{+}$ sont des opérateurs de rang 1, 
satisfaisant
\begin{align*}
E_{+}&=
\frac{1-\chi_{-}}{D_{+}}e_{+}+\MO(h^{\infty}),\\
E_{-}&=
\frac{1}{D_{-}}\langle \bullet, 
\psi_{-}e_{-}\rangle+\MO(h^{\infty}),
\end{align*}
où $D_{\pm}=1+\MO(h^{\infty}),$ et $e_{+},$ $e_{-}$ 
sont normalisés, et
vérifient
\[e_{+}(x,\alpha;h)=
\frac{(\frac{1}{2i}\{\xi+g,\xi+\bar{g}\}(\rho_{+}))^{\frac{1}{4}}}{(\pi h)^{1/4}}
(1+\MO(\tilde{h}))e^{-\frac{i}{h}\int_{x_{+}}^{x}g(y,\alpha)dy},
\]
où $(\xi+g)(\rho_{+})=0$ avec $(\frac{1}{2i}\{\xi+g,\xi+\bar{g}\}(\rho_{+}))>0.$
Nous avons une estimation similaire pour $e_{-}.$ 
Donc 
\begin{align*}
\|(hD_{x}-\widetilde{g}(x,\alpha;h))^{-1}\|
\sim&  \frac{\sqrt{\pi}\exp(\frac{1}{h}\im \ell_{0}(\alpha))}{h^{1/2}
(\frac{1}{2i}\{\xi+g,\xi+\bar{g}\}(\rho_{+}))^{\frac{1}{4}}
(\frac{1}{2i}\{\xi+\bar{g},\xi+g\}(\rho_{-}))^{\frac{1}{4}}}\\
&\times(1+\mathcal{O}( \tilde{h}))
+\mathcal{O}(\frac{1}{\sqrt{h}\,\alpha^{1/4}}),
\end{align*}
où $\tilde{h}:=h/\alpha^{3/2},$ et $\ell_{0}$ (est une intégrale d'action)
vérifie 
\begin{equation}
\ell_{0}:=\im \int_{x_{-}}^{x_{+}}g(x,\alpha) dx.
\end{equation}
\end{theo} 

\section{Cas général en dimension 1}
Soit $p\in S(\R^{2},m),$ indépendant de $h.$ On choisit un 
point $z_{0}$ appartenant 
au bord de l'image du symbole $\Sigma(p)=p(T^{\ast}\R),$
qui est fermé grâce à la condition d'ellipticité  sur $p$ et la croissance à 
l'infini de la fonction d'ordre $m$. 
Nous supposons que 
$p$ est de type principal en $z_{0},$ c'est-à-dire 
\[p(\rho)-z_{0}=0\Rightarrow dp(\rho)\ne0,\]
et que $p$ est d'ordre 2 en  $z_{0}.$

Pour simplifier, nous supposons 
que $p^{-1}(z_{0})$ ne contient qu'un point $\rho_{0}$. 
\\

En nous appuyant sur Hager (\cite{Hager3}):
pour $z_{0}\in \partial \Sigma$ nous avons $p^{-1}(z_{0})=\{\rho_{0}\},$
avec $\{p,\bar{p}\}(\rho_{0})=0.$ Soit $H_{\frac{1}{2i}\{p,\bar{p}\}} $
le vecteur réel  tangent à l'ensemble $\{\rho|\, \frac{1}{2i}\{p,\bar{p}\}(\rho)=0\}.$
Donc 
\begin{equation}\label{gg.1}
(-T_{1},T_{0})\ni s\mapsto \rho(s):=\exp(sH_{\frac{1}{2i}\{p,\bar{p}\}} )(\rho_{0})
\subset\{
\frac{1}{2i}\{p,\bar{p}\}=0\}
\end{equation}
est une courbe orientée,
et $\frac{1}{2i}\{p,\bar{p}\}$ est positif à la gauche de cette courbe, et 
négatif à sa droite.
Concernant la ligne (\ref{gg.1}),
nous avons évidemment $\rho(0)=\rho_{0}$,
et le sous-ensemble de $\C,$
$p\left(
\{\frac{1}{2i}\{p,\bar{p}\}=0\}
\right),$
coïncide avec $\partial\Sigma$ près de $z_{0}$.
\\

$p(\rho(s))$ est aussi une courbe orientée, puisque
\[
\frac{\partial}{\partial s} (p(\rho(s)))=
(H_{\frac{1}{2i}\{p,\bar{p}\}}p)(\rho(s))\\
=
-\{p,\frac{1}{2i}\{p,\bar{p}\}\}(\rho(s))\ne 0.
\]
De plus $d\im p\wedge d\re p =\frac{1}{2i}\{p,\bar{p}\}d\xi\wedge dx,$ donc 
$\rho\to p(\rho)$ préserve l'orientation à gauche de $\rho(s),$ et la renverse à droite.
Ce qui veut dire,  qu'au voisinage de $z_{0},$ 
si $\rho(s)$ est parcouru dans le sens positif,
on voit que  $\Sigma$ se trouve à 
gauche de la courbe $p(\rho(s))\in \partial\Sigma.$

On note le chemin $p(\rho(s))$ par $\gamma: ]-T_{1},T_{0}[\to \C.$ 
Un vecteur tangent à $\partial \Sigma$ est 
donné par $\dot{\gamma}(s),$
et nous savons que
\begin{equation}
\dot{\gamma}(s)=
-\{p,\frac{1}{2i}\{p,\bar{p}\}\}(\rho(s)).
\end{equation}
Il existe un voisinage de $z_{0}$ où tout point $z$ appartenant 
à $\Sigma(p)$ s'écrit sous la forme
\begin{equation}\label{gg.3}
z=\gamma(s)+\alpha n(s),\quad 0\le\alpha<\frac{1}{C_{0}},
\quad -T_{1}<s<T_{0},
\quad z_{0}=\gamma(0),\\
\end{equation}
pour une constante $C_{0}>0$ donnée, où $n(s)$ est le vecteur unitaire normal 
à la courbe $\gamma(s),$ orienté vers l'intérieur de $\Sigma(p).$ 
$n(s)$ vérifie les relations $n(s)=i\frac{\dot{\gamma}(s)}{|\dot{\gamma}(s)|},$
et $|n(s)|=1.$ Par la suite le vecteur unitaire tangent à $\gamma(s)$ sera 
noté $u(s),$ $u(s)=\frac{\dot{\gamma}(s)}{|\dot{\gamma}(s)|}$ et $|u(s)|=1.$

Puisque $p$ est de type fini d'ordre 2 au point $z_{0}$,
nous pouvons alors supposer
sans perte de généralité que 
$p'_{\xi}(\rho_{0})\ne 0$ et que
$p'_{x}(\rho_{0})= 0.$
Nous avons dans de nouvelles coordonnées 
$(x,\xi)$ centrées en $\rho(s)$ - changement de coordonnées symplectiques
 dépendant explicitement de $s$, 
conjugaison par un opérateur intégral de Fourier de plus le gradient de $p$ en 
$(0,0)$ est un multiple d'un vecteur réel -
\begin{align*} 
p(\rho,s)-\gamma(s)&=p(\rho,s)-p(\rho(s))\\
&=u(s)(\frac{p'_{\xi}(0,0,s)}{u(s)}\xi+\omega(x,\xi,s)),
\mbox{ avec } \omega=\MO(x^{2}+\xi^{2}),
\end{align*}
où $\omega$ est $C^{\infty}$ en $s.$ N'oublions pas que $p$ va dépendre de 
$s$ car les nouvelles coordonnées dépendent de $s$. Aussi $p'_{\xi}(0,0,s)/u(s)$ est réel.
Puis, grâce au théorème de factorisation de Malgrange pour des fonctions 
$C^{\infty}$,  il existe un ouvert $V$ de $\R^{2}\times\R_{+}\times]-T_{1},T_{0}[$
contenant le point $(0,0,0,0)$ tel que dans $V,$ nous avons la factorisation 
\begin{align}
p(\rho,s)-z&=p(\rho,s)-\gamma(s)-i\alpha u(s)\nonumber\\
&=u(s)(\frac{p'_{\xi}(0,0,s)}{u(s)}\xi+\omega(x,\xi,s)-i\alpha)\nonumber\\
&=u(s)q(x,\xi,\alpha,s)\,(\xi+g(x,\alpha,s)),\quad
q(x,\xi,\alpha,s)\ne 0,
\end{align}
où $g$ et $q$ sont des fonctions $C^{\infty},$ $g$ s'annulant au point $(0,0,0).$
De plus pour $x=\xi=0,$ $\alpha=0$ et $s$ non nul, nous avons 
$p(0,s)-z=p(0)-\gamma(s)=0,$ impliquant $g(0,0,s)=0.$ 
Par commodité, nous supprimons la variable $s$.
Pour $x=\xi=0,$ $\alpha=0,$ nous avons alors
\begin{align}
& \quad g(0,0)=0, \label{gg.4}\\
& \quad q(0,0,0)= \frac{p'_{\xi}(0,0)}{u(s)}, \\
& \quad  g'_{\alpha}(0,0)=-i \frac{u(s)}{p'_{\xi}(0,0)},\label{gg.4bis} \\
& \quad  g'_{x}(0,0)=0. 
\end{align}
Comme $p'_{x}(\rho_{0})$ est nul, le rapport $\frac{p'_{\xi}(0,0)}{u(s)}$ 
est réel et non nul, et sera noté $a$; $a:=\frac{p'_{\xi}(0,0)}{u(s)}=q(0,0,0).$ 
Puisque
\[\omega''_{xx}(0,0)=q''_{xx}(0,0,0)\,.\,0+2\,q'_{x}(0,0,0)\,.\,0+
q(0,0,0)g''_{xx}(0,0),\]
nous avons aussi,
\begin{equation}
 a g''_{xx}(0,0)=\omega''_{xx}(0,0).
\end{equation}
Nous avons alors - rappelons que c'est pour $\alpha=0$ -
\begin{align*}
\frac{1}{2i}\{p,\bar{p}\}(x,\xi)&=
\frac{|u(s)|^{2}}{2i}\{a\xi+\omega,a\xi+\bar{\omega}\}(x,\xi)\\
&=\frac{1}{2i}(a\overline{\omega'_{x}}-a\omega'_{x})\\
&=-a \im \omega'_{x}(x,\xi),
\end{align*}
ce qui implique
\begin{align}
\{p,\frac{1}{2i}\{p,\bar{p}\}\}(0,0)&=
u(s)
\{a\xi+\omega,
-\im a\omega'_{x}\}(0,0)\nonumber\\
&=-u(s)a^{2}\im \omega''_{xx}(0,0)\nonumber\\
&=-u(s)a^{2}a\im g''_{xx}(0,0).
\label{gg.5}
\end{align}
On déduit alors que
\begin{equation}
\im g''_{xx}(0,0)\ne 0,
\end{equation}
et que 
\begin{align}
\im g''_{xx}(0,0)=\frac{\dot{\gamma}(s)}{u(s)a^{2}a}
&=\frac{\dot{\gamma}(s)^{3}}{|\dot{\gamma}(s)|^{3}}\frac{|\dot{\gamma}(s)|}{p'_{\xi}(0,0)^{3}}\nonumber\\
&=\frac{\{p,\frac{1}{2i}\{p,\bar{p}\}\}(0,0)^{2}}{|\{p,\frac{1}{2i}\{p,\bar{p}\}\}(0,0)|^{2}}
\frac{\{p,\frac{1}{2i}\{p,\bar{p}\}\}(0,0)}{p'_{\xi}(0,0)^{3}}.\label{gg.5bis}
\end{align}

De l'équation (\ref{g.0}) qui donne le comportement de $x_{\pm}(\alpha)$  dans le cas modèle $hD_{x}+g(x,\alpha)$ avec l'hypothèse 
que  $g'_{\alpha}(0,0)=-i$, nous obtenons
qu'il existe un voisinage de $z_{0}$ pour lequel nous avons pour
tout complexe $z$ appartenant à $\Sigma(p)$ 
 \begin{equation}
z\notin \partial\Sigma(p),\; p^{-1}(z)=\{\rho_{+}(z),\rho_{-}(z)\},\quad
\rho:=(x,\xi),
\end{equation}
où 
\[
x_{\pm}(z)=x_{0}\mp \epsilon\-(\alpha|g'_{\alpha}(0,0,s)|)^{1/2}
\left(\frac{2}{|\im g''_{xx}(0,0,s)|}\right)^{1/2}
+\MO(\alpha),
\]
$\epsilon$ est le signe de $g'_{\alpha}(0,0,s)\im g''_{xx}(0,0,s)$. En rempla\c cant
$g'_{\alpha}(0,0,s)$ et $\im g''_{xx}(0,0,s)$ par leur expression respective  (\ref{gg.4bis}) et 
(\ref{gg.5bis}),
nous obtenons
\begin{align}
x_{\pm}(z)&=
x_{0}\mp \epsilon\,\alpha^{1/2}|p'_{\xi}(\rho(s))|
\left(\frac{2}{|\{p,\frac{1}{2i}\{p,\bar{p}\}\}(\rho(s))|}
\right)^{1/2}
+\MO(\alpha),\label{gg.6}\\
z&=\gamma(s)+\alpha n(s), \alpha>0,\nonumber
\end{align}
où $\epsilon$ est le signe de $\{p,\frac{1}{2i}\{p,\bar{p}\}\}(\rho(s))/p'_{\xi}(\rho(s))$.
Nous avons que  
\begin{equation}
\pm\frac{1}{2i}\{p,\bar{p}\}(\rho_{\pm})> 0.
\end{equation}

On prolonge $g$ sur $\R$ de la même manière que dans le cas modèle à paramètre.\\

On distingue, maintenant, deux cas :
\begin{itemize}
\item[$\bullet$] 1er cas,  $p^{-1}(z_{0})$ ne contient qu'un point $\{\rho_{0}\}$,
\item[$\bullet$] 2ème cas,  $p^{-1}(z_{0})$ ne contient que deux points $\{\rho_{1},\rho_{2}\}$.
\end{itemize}
\subsection{Estimation de résolvante: cas 1} 
%

\begin{prop}
Il existe deux fonctions
\begin{align}
\widetilde{q}&\sim \sum_{k\ge 0} \widetilde{q}_{k}h^{k}\in S_{cl}(\mathrm{Vois}(0,0),1),\\
\widetilde{g}&\sim \sum_{k\ge 0} \widetilde{g}_{k}h^{k},\quad g_{k}\in C^{\infty}(\pi_{x}(\mathrm{Vois}(0,0))),
\end{align}
telles que - nos coordonnées $(x,\xi)$ sont centrées en $\rho(s)$ - 
\begin{equation}\label{gg.7}
p(x,\xi)-z\sim
\widetilde{q}(x,\xi,\alpha,s;h)\#(\xi+\widetilde{g}(x,\alpha,s;h))
\#\widetilde{q}(x,\xi,\alpha,s;h),
\end{equation}
dans $S(\mathrm{Vois}(0,0),m)$ uniformément en $z$ ($z$ dépend de $s$ et $\alpha$), 
où $z$ appartient à un voisinage de $z_{0}.$ Si l'on tient compte de la factorisation
$(p-z)=n(s)q(x,\xi)(\xi+g(x)),$ nous pouvons choisir  $\widetilde{q}_{0}$ égal à  $(n(s)q(x,\xi))^{1/2}$
impliquant que $\widetilde{g}_{1}=0.$
\end{prop}
Pour la suite, par souci de clarté, nous supprimons la dépendance en $\alpha$ et $s$.
\noindent\textbf{Preuve.} Avec $p(\rho(z))-z=0,$ $p'_{\xi}(\rho(z))\neq 0$ nous obtenons
un voisinage $V$ de $\rho(z_{0})$ des fonctions lisses $\widetilde{q}_{0}$  et  $\widetilde{g}_{0}$
définies respectivement sur $V$ et $\pi_{x}(V)$ telles que 
\[(p-z)=(n(s)q(x,\xi))^{1/2}(\xi+\widetilde{g}_{0}(x))(n(s)q(x,\xi))^{1/2}\]
avec $\widetilde{q}_{0}(x(z_{0}),\xi(z_{0}))\ne 0,$ $\widetilde{g}_{0}(x(z_{0})=-\xi(z_{0})$
On peut rajouter $z$ aux variables et nous avons toujours une dépendance $C^{\infty}$ de $z$, 
et les équations ci-dessus restent valables dans un voisinage de $z_{0}.$
Nous regroupons par ordre de $h$ les termes de la forme de composition asymptotique:
\begin{align*}
0=&\widetilde{q}_{0}^{2}(x,\xi) \widetilde{g}_{N}(x)+2\widetilde{q}_{0}(x,\xi)
\widetilde{q}_{N}(x,\xi)(\xi+\widetilde{g}_{0}(x))\\
&+\widetilde{G}_{N}(\widetilde{q}_{0},\ldots,\widetilde{q}_{N-1},
\widetilde{g}_{0},\dots \widetilde{g}_{N-1},x,\xi,z).
\end{align*}
Etant donné que $\widetilde{q}_{0}$ est non nul, nous avons une équation de la forme
\[
G_{N}=2\frac{\widetilde{q}_{N}}{\widetilde{q}_{0}}(\xi+\widetilde{q}_{0}(x))
+\widetilde{g}_{N}(x)
\]
où $G_{N}$ ne dépend que de $\widetilde{q}_{0},\ldots,\widetilde{q}_{N-1},
\widetilde{g}_{0},\dots \widetilde{g}_{N-1}$. Il est possible de déterminer 
$\widetilde{g}_{N}$ et $\widetilde{q}_{N}$ inductivement (voir \cite{Hager2}), puisque le
théorème de Malgrange nous prouve l'existence d'un voisinage de $\rho(z_{0})$ et les 
fonctions $\widetilde{g}_{N}$ et $\widetilde{q}_{N}$ avec les propriétés demandées. 
En itérant, nous obtenons les séries formelles $\widetilde{g}$ et $\widetilde{q}$. 

Pour le dernier point de la proposition, il suffit de remarquer que le terme en $h$
d'une composition de Weyl symétrique s'annule.
\hfill $\square$ \medskip


Nos variables $\rho=(x,\xi)$ sont centrées en $\rho(s):=(x(s),\xi(s))$; si nous introduisons 
l'opérateur unitaire $T_{s}$ par 
\[ T_{s}u(y):=e^{iy\xi(s)/h}u(y-x(s)),\]
nous trouvons que l'opérateur $\widetilde{P}:=T_{s}PT^{-1}_{s}$ satisfait
\[\widetilde{P}=\mathrm{Op}^{w}_{h}(\widetilde{p})\mbox{ où } \widetilde{p}(\rho)=p(\rho-\rho(s)).\]

Il existe alors un symbole $\widetilde{q}$ 
pour les variables $(x,\xi)$, $\alpha$ et $s,$ vérifiant 
$\widetilde{q}=\chi(x,\xi)q(x,\xi,\alpha,s)^{-\frac{1}{2}}+\MO(h),$
où $\chi=C^{\infty}_{0}(\pi_{x,\xi}(V))$ est une troncature à support compact 
indépendant de $z$ (donc de $\alpha$ et $s$)
valant $1$ dans un petit voisinage $V$ de $(0,0),$ tel que
\begin{align}
(\widetilde{q})^wT_{s}(P-z)T_{s}^{-1}(\widetilde{q})^w=n(s)
& (hD_{x}+g(x,\alpha,s)+h^{2}g_{\ast}(x,\alpha,s;h)) \nonumber\\
&\mbox{ microlocalement près de } (0,0),
\end{align}
\textit{(remarque de J. Sj\"ostrand, pour y arriver, il semblerait qu'il faut aussi une 
rotation symplectique)} 
uniformément pour $z$ (qui est une fonction de $\alpha$ et $s,$ voir (\ref{gg.3})) 
dans un voisinage donné de $z_{0},$
\begin{align*}
g_{\ast}(x,\alpha,s;h)\sim g^{0}_{\ast}(x,\alpha,s)+hg^{1}_{\ast}(x,\alpha,s)+\ldots
\mbox{ dans } C^{\infty}_{0},\\
\mbox{ uniformément en } \alpha \mbox{ et } s.
\end{align*}

Par souci de clarté, la variable $s$ est supprimée. On notera 
\[\widetilde{g}(x,\alpha;h):=g(x,\alpha)+h^{2}g_{\ast}(x,\alpha;h).\]
%
%


Nous pouvons trouver une partition de l'unité de $\R^{2},$
$\chi_{j}\in S(\R^{2},1),$ 
$ j =0,1$  telle que (voir Faure-Sj\"ostrand par exemple)
\[ 
\chi_{j}\ge 0,\;
\chi_{0}^{2}+\chi_{1}^{2}=1,\mbox{ avec }
\chi_{0}\in C^{\infty}_{0}(\mbox{Vois}(\rho_{0})),\,
 \chi_{1}\in C^{\infty}_{b}.
\]
Nous avons
\[
\|u\|^{2}=\|\chi_{1}u\|^{2}+\|\chi_{0}u\|^{2}.
\]

En utilisant l'ellipticité de $(P-z)^{-1}$ on peut montrer que pour 
$|\alpha |=\MO (h\ln(\frac{1}{h})^{2/3})$
\begin{equation}
\|(P-z)^{-1}v\|^{2}\le M_{0}^{2}\|\chi_{0}v\|^{2}+\MO(1)\|\chi_{1}v\|^{2}+\MO(h)\|v\|^{2},
\end{equation}
où $M_{0}=\|(P-z)^{-1}\chi_{0}\|.$ Ce qui veut dire que si $M:=\max(M_{0},\MO(1))$ alors
\[\|(P-z)^{-1}\|\le (M^{2}+\MO(h))^{\frac{1}{2}}.\]

Nous souhaitons maintenant connaître
le comportement de $M_{0}$,
\[
M_{0}=\| 
(n(s))^{-1}T_{s}^{-1}\widetilde{q}^w(hD_{x}+\widetilde{g}(x,\alpha;h))^{-1}
\widetilde{q}^wT_{s}\chi_{0}\|.
\]
 Du théorème \ref{g23} et du fait que 
$E_{-+}$ est scalaire et $T_{s}$ est unitaire, nous avons que 
\begin{align*}
M_{0}= |E_{-+}^{-1}|\times \| \widetilde{q}^w
E_{+}E_{-} \widetilde{q}^w\|
+\mathcal{O}(\frac{1}{\sqrt{h}\,\alpha^{1/4}}),
\end{align*}
et nous obtenons une expression 
pour $\widetilde{q}^wE_{+}E_{-} \widetilde{q}^w,$
en fonction des quasimodes:
\begin{equation}
\widetilde{q}^w
E_{+}E_{-}  \widetilde{q}^w
=\frac{1-\chi_{-}}{D_{-}D_{+}}(\widetilde{q}^{w}e_{+})\langle \bullet, 
((\widetilde{q})^\ast) ^w\psi_{-}e_{-}\rangle+\MO(h^{\infty}).
\end{equation}
Nous avons besoin  d'un théorème de Melin, Sj\"ostrand 
sur l'action d'un opéra\-teur pseudodifférentiel sur une fonction BKW
avec phase complexe. Dès lors, puisque $e_{+}$ et $e_{-}$ sont des 
fonctions BKW normalisées, et 
microlocalisées près de $\rho_{\pm},$ nous avons le résultat 
suivant:
\begin{lemme}
Nous avons la formule 
 \begin{equation}
\|\widetilde{q}^w
E_{+}E_{-} \widetilde{q}^w\|
\sim \frac{1}{|q(\rho_{+})|^{\frac{1}{2}}
|q(\rho_{-})|^{\frac{1}{2}}}(1+\MO(\tilde{h})).
\end{equation}
\end{lemme}
Nous utilisons maintenant l'expression (\ref{g.60})
de $E_{-+}.$ 
Puisque  $(p-z)=n(s)q(\xi+g)$, avec $|n(s)|=1$ nous obtenons que
\[\frac{1}{|q(\rho_{+})|^{\frac{1}{2}}}\times
(\frac{1}{2i}\{\xi+g,\xi+\bar{g}\}(\rho_{+}))^{\frac{1}{4}}
=
(\frac{1}{2i}\{p,\bar{p}\}(\rho_{+}))^{\frac{1}{4}},
\]
et similairement pour $\rho_{-}$. 
Nous avons donc,
pour tout  $h\ll \alpha^{3/2}\ll 1,$ assez petit  
\begin{equation}\label{gg.10}
M_{0}=\frac{\sqrt{\pi}\exp(\frac{1}{h}\im \ell_{0}(\alpha))}{h^{1/2}
(\frac{1}{2i}\{p,\bar{p}\}(\rho_{+}))^{\frac{1}{4}}
(\frac{1}{2i}\{\bar{p},p\}(\rho_{-}))^{\frac{1}{4}}}
(1+\mathcal{O}( \tilde{h}))
+\mathcal{O}(\frac{1}{\sqrt{h}\,\alpha^{1/4}}),
\end{equation} 
où $\tilde{h}:=h/\alpha^{3/2},$ et $\ell_{0}$ (qui est une intégrale d'action)
vérifie 
\begin{equation}
\ell_{0}(\alpha):=\im \int_{x_{-}}^{x_{+}}\im g(x,\alpha,s) dx
\end{equation}
(remarquons l'ajout de $s$ qui est important pour la suite). 
$\ell_{0}(\alpha)=\ell_{0}(z)$ s'écrit aussi sous la forme 
\begin{align}
\ell_{0}(z)
&=\im \int_{x_{-}(z)}^{x_{+}(z)}\varphi(x,z) dx, \quad p(x,\varphi(x,z))-z=0,\nonumber\\
&=\im \int_{\gamma\subset p^{-1}(z)}\xi dx,
(\gamma \mbox{ relie } \rho_{-}\mbox{ à }\rho_{+}),\label{gg.11}\\
&=\im \int_{\gamma\subset p^{-1}(z)}x d\xi\nonumber .
\end{align}
Pour que $M_{0}$ ait que croissance tempérée en $h^{-1}$, il faut que 
\[(\alpha^{3/2}\asymp)\quad
\im \int^{x_{+}}_{x_{-}}g(y,z)dy\le Ch\ln \frac{1}{h},\]
pour un $C>0$ quelconque. Soit si 
\begin{equation}
\alpha\le \MO(1) \left( h\ln\frac{1}{h} \right)^{2/3}.
\end{equation}

On déduit alors le résultat suivant :
\begin{theo}
Soit $P$ un opérateur $h$-pseudodifférentiel de symbole indépen\-dant de $h.$
On suppose que $(p,z_{0})$ vérifie les conditions (\ref{int.9}), (\ref{int.10}),
(\ref{int.13}).
Un vecteur tangent à $\partial \Sigma(p)$ dans un voisinage de $z_{0}$ est 
donné par $\dot{\gamma}(s).$
Donc tout point $z$ de $\Sigma(p)$ appartenant à un voisinage de $z_{0}$ peut s'écrire
sous la forme
\begin{equation}
z=\gamma(s)+\alpha n(s),\quad -T_{1}<s<T_{0},\quad  \alpha\ge 0,
\quad z_{0}=\gamma(0),
\end{equation}
où $n(s)=i\dot{\gamma}(s)/|\dot{\gamma}(s)|.$ Il existe un voisinage de $z_{0},$ 
pour lequel tout point de $z$ à l'intérieur
de $\Sigma(p)$ vérifie
\begin{equation}\label{gg.13}
z\notin \partial\Sigma,\; p^{-1}(z)=\{\rho_{+}(z),\rho_{-}(z)\},
\quad \rho_{\pm}=(x_{\pm},\xi_{\pm}),
\end{equation}
avec
\[\pm\frac{1}{2i}\{p,\bar{p}\}(\rho_{\pm})>0.\]
Il existe une constante $T_{\ast}$ ($<T_{0},T_{1}$) telle que pour toutes constantes 
$C_{0},C_{1}>0,$ il existe une constante $C_{2}>0$ telle que la résolvante 
$(P-z)^{-1}$ est bien définie pour 
\begin{equation}
 |s|<T_{\ast},\;
\frac{ h^{2/3}}{C_{0}}\le \alpha\le C_{1}(h\ln \frac{1}{h})^{2/3},\;
 h<\frac{1}{C_{2}},
 \end{equation}
 et satisfait l'estimation
\begin{align*}
\|(P-z)^{-1}\|
\sim&  \frac{\sqrt{\pi}\exp(\frac{1}{h}\im \ell_{0}(z))}{h^{1/2}
(\frac{1}{2i}\{p,\bar{p}\}(\rho_{+}))^{\frac{1}{4}}
(\frac{1}{2i}\{\bar{p},p\}(\rho_{-}))^{\frac{1}{4}}}
\times(1+\mathcal{O}( \tilde{h}))
+\mathcal{O}(\frac{1}{\sqrt{h}\,\alpha^{1/4}}),
\end{align*}
où $\tilde{h}=h/\alpha^{3/2},$ et $\ell_{0}$ (qui est une intégrale d'action)
vérifie 
\begin{equation}
\ell_{0}(z):=\im \int_{\gamma\subset p^{-1}(z)}\xi dx,
(\gamma \mbox{ relie } \rho_{-}\mbox{ à }\rho_{+}).
\end{equation}
\end{theo}

Dans le cas des opérateurs à coefficients analytiques, on peut pousser 
$\alpha$ jusqu'à une petite constante.

\subsection{Estimation de résolvante: cas 2}
\begin{theo}
Soit $P$ un opérateur $h$-pseudodifférentiel de symbole indépen\-dant de $h.$
On demande que pour tout $z$ dans l'intérieur de $\Sigma,$ $p^{-1}(z)$ soit 
fini et que $\{p,\bar{p}\}(\rho)\ne 0$ pour $\rho\in p^{-1}(z).$

On suppose que p vérifie les conditions (\ref{int.9}), (\ref{int.10}), et 
que l'hypothèse  (\ref{int.13}) soit remplacée par 
\begin{align}\label{egg.10}
p^{-1}(z_{0})=&\{\rho_{1},\rho_{2}\},\,
\rho_{j}=(x_{j},\xi_{j}),\nonumber\\
&\{p,\{p,\bar{p}\}\}(\rho_{1})\ne 0,\nonumber\\
&\{p,\{p,\bar{p}\}\}(\rho_{2})\ne 0.
\end{align}

$\Sigma(p)$ sera l'union de deux domaines $\Sigma_{1},$ $\Sigma_{2}$ à bord $C^{\infty}$
près de $z_{0}$ dont les bords peuvent ne pas coïncider, $\Sigma_{j}=p(\mbox{vois}(\rho_{j})).$

Soit $H_{\frac{1}{2i}\{p,\bar{p}\}} $
le vecteur réel  tangent à l'ensemble $\{\rho|\, \frac{1}{2i}\{p,\bar{p}\}(\rho)=0\}.$
Donc 
\[(-T_{1},T_{0})\ni s\mapsto \rho_{j}(s):=\exp(sH_{\frac{1}{2i}\{p,\bar{p}\}} )(\rho_{j})\subset\{
\frac{1}{2i}\{p,\bar{p}\}=0\}\]
est une courbe orientée. 
Près de $z_{0},$ pour tout $z_{\ast}\in\partial\Sigma(p),$ nous avons
$p^{-1}(z)=\{\rho_{1}(s_{\ast}),\rho_{2}(t_{\ast})\}.$
$p(\rho_{j}(s))$ est aussi une courbe orientée 
notée $\gamma_{j}(s):=p(\rho_{j}(s)).$
Le vecteur tangent à $\partial \Sigma$ dans un voisinage 
de $z_{0}$ est donné par 
$\dot{\gamma}_{1}(s)$ ou par $\dot{\gamma}_{2}(t).$
Il existe donc un voisinage de $z_{0}$ tel que tout point $z$ de $\Sigma(p)$ 
se met sous la forme :
\begin{align}
z&=\gamma_{1}(s)+\alpha_{1} n_{1}(s),\quad-T_{1}<s<T_{0},\quad \alpha_{1}\ge 0,
\quad z_{0}=\gamma_{1}(0),\\
z&=\gamma_{2}(s)+\alpha_{2} n_{2}(s)\quad,-T_{1}<s<T_{0},\quad \alpha_{2}\ge 0,
\quad z_{0}=\gamma_{2}(0),
\end{align}
$n_{1}=i\dot{\gamma}_{2}(s)/|\dot{\gamma}_{1}(s)|$ et 
$n_{2}=i\dot{\gamma}_{2}(s)/|\dot{\gamma}_{1}(s)|$.
Alors pour tout point $z$ de $\Sigma$ dans un voisinage de $z_{0},$
nous avons 
\begin{equation*}
z\notin \partial\Sigma,\; p^{-1}(z)=\{\rho_{+}^{1}(z),\rho_{-}^{1}
(z), \rho_{+}^{2}(z),\rho_{-}^{2}(z)\},
\end{equation*}
où
\[\pm\frac{1}{2i}\{p,\bar{p}\}(\rho_{\pm}^{j})>0.\]

Il existe une constante $T_{\ast}$ ($<T_{0},T_{1}$) telle que pour toutes constantes 
$C_{0},C_{1}>0$ il existe une constante $C_{2}>0$ telle que la résolvante 
$(P-z)^{-1}$ est bien définie pour 
\begin{equation}
 |s|<T_{\ast},\;
\frac{ h^{2/3}}{C_{0}}\le \alpha\le C_{1}(h\ln \frac{1}{h})^{2/3},\;
 h<\frac{1}{C_{2}},\; \alpha=\max(\alpha_{1},\alpha_{2}),
 \end{equation}
 et satisfait l'estimation pour $z\in\Sigma_{1}\cap\Sigma_{2}:$
\begin{align*}
\|(P-z)^{-1}\| \sim \sup_{i=1,2} \bigg(
&  \frac{\sqrt{\pi}\exp(\frac{1}{h}\im \ell_{0}^{i}(z))}{h^{1/2}
(\frac{1}{2i}\{p,\bar{p}\}(\rho_{+}^{i}))^{\frac{1}{4}}
(\frac{1}{2i}\{\bar{p},p\}(\rho_{-}^{i}))^{\frac{1}{4}}}\bigg)
\times(1+\mathcal{O}( \tilde{h}))\\
&+\mathcal{O}(\frac{1}{\sqrt{h}\,\alpha^{1/4}}),
\end{align*}
pour $z\in\Sigma_{1}\cap\Sigma^{c}_{2}:$
\begin{align*}
\|(P-z)^{-1}\| \sim 
&  \frac{\sqrt{\pi}\exp(\frac{1}{h}\im \ell_{0}^{1}(z))}{h^{1/2}
(\frac{1}{2i}\{p,\bar{p}\}(\rho_{+}^{1}))^{\frac{1}{4}}
(\frac{1}{2i}\{\bar{p},p\}(\rho_{-}^{1}))^{\frac{1}{4}}}
\times(1+\mathcal{O}( \tilde{h}))\\
&+\mathcal{O}(\frac{1}{\sqrt{h}\,\alpha^{1/4}}),
\end{align*}
pour $z\in\Sigma_{1}^{c}\cap\Sigma_{2}:$
\begin{align*}
\|(P-z)^{-1}\| \sim 
&  \frac{\sqrt{\pi}\exp(\frac{1}{h}\im \ell_{0}^{2}(z))}{h^{1/2}
(\frac{1}{2i}\{p,\bar{p}\}(\rho_{+}^{2}))^{\frac{1}{4}}
(\frac{1}{2i}\{\bar{p},p\}(\rho_{-}^{2}))^{\frac{1}{4}}}
\times(1+\mathcal{O}( \tilde{h}))\\
&+\mathcal{O}(\frac{1}{\sqrt{h}\,\alpha^{1/4}}),
\end{align*}
où $\tilde{h}=h/\alpha^{3/2},$ et $\ell_{0}^{i}$ (qui est une intégrale d'action)
vérifie 
\begin{equation}
\ell_{0}^{i}(z):=\im \int_{\gamma\subset p^{-1}(z)}\xi dx,
(\gamma \mbox{ relie } \rho_{-}^{i}\mbox{ à }\rho_{+}^{i}).
\end{equation}
\end{theo}
\noindent\textbf{Preuve.} Fort du cas 1, nous pouvons donner les grandes lignes de la preuve.
Comme précédemment notre raisonnement  donne donc deux opérateurs 
$T_{s}$ et $T_{s}^{-1}$ tels que 
 \begin{align}
 \widetilde{q}^wT_{s}(P-z)T_{s}^{-1}\widetilde{q}^w=
\nonumber&
n_{1}(s) (hD_{x}+\widetilde{g}_{1}(x,\alpha;h))\\
&\mbox{ microlocalement près de } (0,0),
\end{align}
\[\widetilde{g}_{1}=g_{1}(x,\alpha)+h^{2}g_{\ast,1}(x,\alpha;h),\]
De même pour $\rho_{2}.$ Soient $V_{1},V_{2}$ 
des voisinages arbitrairement petits près de $\rho_{1},\rho_{2}$. Il existe 
$\theta_{j}\in C^{\infty}_{0}(V_{j})$ tels que si
\[p_{0}=p+\sum_{j}\theta_{j}\]
alors $p_{0}^{-1}(z_{0})$ est vide. Soit $p_{j}=p+\sum_{k,k\ne j}\theta_{k}.$
Donc $p^{-1}_{j}(z_{0})=\{\rho_{j}\}.$ 
Soit $P_{0}=p_{0}^{w}(x,hD_{x})$ et $P_{j}=p_{j}^{w}(x,hD_{x})$.
$P_{0}$ est uniformément elliptique. On suppose que $\theta_{0}=1$
près de $\rho_{0}.$ Nous écrivons
\[
R_{0}(z):= \sum_{j}(P_{j}-z)^{-1}\theta_{j}+(P_{0}-z)^{-1}(1-\sum_{j}\theta_{j})
\]

Nous avons 
\begin{align*}
(P-z)R_{0}(z)=&1+\sum_{j}(P-P_{j})(P_{j}-z)^{-1}\theta_{j}\\
&+(P+P_{0})(P_{0}-z)^{-1}(1-\sum_{j}\theta_{j}),
\end{align*}
où $(P-P_{j})(P_{j}-z)^{-1}\theta_{j}=\MO(h^{\infty}),$ (voir \cite{Hager2}) et aussi 
pour le dernier terme.
Donc $(P-z)=R_{0}(z)(1+K)^{-1}$ où $K=\MO(h^{\infty}).$

Nous pouvons trouver une partition de l'unité de $\R^{2},$
$\chi_{j}\in S(\R^{2},1),$ 
$ j =0,1,2$  telle que 
\[ 
\chi_{j}\ge 0,\;
\chi_{0}^{2}+\ldots+ \chi_{2}^{2}=1,\mbox{ avec }
\chi_{j}\prec \theta_{j},\,
 \chi_{0}\in C^{\infty}_{b}.
\]
Nous avons
\[
\|u\|^{2}=\sum_{j=0}^{2}\|\chi_{j}u\|^{2}.
\]
On peut alors montrer que
%
\begin{equation}
\|(P-z)^{-1}v\|^{2}\le \sum_{j=1}^{2} M_{j}^{2}\|\chi_{j}v\|^{2}+\MO(1)\|\chi_{0}v\|^{2}+\MO(h)\|v\|^{2},
\end{equation}
où $M_{j}=\|(P_{j}-z)^{-1}\chi_{j}\|.$ Ce qui veut dire que si $M:=\max(M_{j},\MO(1))$ alors
\[\|(p-z)^{-1}\|\le (M^{2}+\MO(h))^{\frac{1}{2}}.\]

Nous souhaitons maintenant connaître
le comportement de $M_{j}$ qui est aussi donné par l'expression,
\begin{align}
M_{1}&=\| 
(n_{1}(s))^{-1}T_{s}^{-1}\widetilde{q}^w_{1}(hD_{x}+\widetilde{g}_{1}(x,\alpha_{1};h))^{-1}
\widetilde{q}^w_{1}T_{s}\chi_{1}\|,\nonumber\\
M_{2}&=\|
(n_{2}(s))^{-1}S_{s}^{-1}\widetilde{q}^w_{2}(hD_{x}+\widetilde{g}_{2}(x,\alpha_{2};h))^{-1}
\widetilde{q}^w_{2}S_{s}\chi_{2}\|,
\end{align}
$S_{t}$ a un comportement similaire à $T_{s}$.

De la section précédente, $M_{j},$ $j=1,2$  satisfait  
\begin{equation}
M_{j}\sim  \frac{\sqrt{\pi}\exp(\frac{1}{h}\im \ell_{0}^{j}(z))}{h^{1/2}
(\frac{1}{2i}\{p,\bar{p}\}(\rho_{+}^{j}))^{\frac{1}{4}}
(\frac{1}{2i}\{\bar{p},p\}(\rho_{-}^{j}))^{\frac{1}{4}}}
(1+\mathcal{O}( \tilde{h}_{j}))
+\mathcal{O}(\frac{1}{\sqrt{h}\,\alpha^{1/4}}),
\end{equation} 
où $\tilde{h}_{j}=h/\alpha_{j}^{3/2},$  et $\ell_{0}^{j}$ 
(qui est une intégrale d'action)
vérifie 
\begin{equation}
\ell_{0}^{j}(z):=\im \int_{\gamma\subset p^{-1}(z)}\xi dx,
(\gamma \mbox{ relie } \rho_{-}^{j}\mbox{ à }\rho_{+}^{j}).
\end{equation}
\hfill $\square$ \medskip

\section{Exemples}
Dans ce qui précède, le paramètre $\alpha$ représentait la distance 
du point $z$ à la frontière du pseudospectre. Dans ce qui suit, pour les exemples, 
la frontière ne sera représentée que par des morceaux de l'axe réel ou imaginaire.
Donc la partie réelle ou imaginaire de $z$ sera notre distance à la frontière, soit notre paramètre 
$\alpha.$ On préfère ici changer et utiliser la lettre $y$ à la place de $\alpha.$ 
\subsection{Cas 1: L'opérateur cubique non-autoadjoint}
On considère l'opérateur non-autoadjoint sur la droite réelle
\[\mathcal{P}=(D_{x})^{2}+ix^{3}.\]

L'image du symbole principal est le demi-plan $\re z\ge 0.$
Le spectre est composé de valeurs propres discrètes situées sur 
l'axe réel positif (voir Trepheten, Embree \cite{TE}).
Le pseudospectre est symétrique
par rapport à l'axe réel. Pour étudier les courbes de niveaux de la résolvante,
il suffit donc de considérer celles situées au-dessus de l'axe réel.
$p$ désigne le symbole de $\mathcal{P}$.\\

Pour $z$ au dessus de l'axe réel et appartenant à $\partial\Sigma(p)=i\R,$
nous avons $p^{-1}(z)=\{ \rho \}=\{(\im z)^{\frac{1}{3}},0\}$ et
\[\{p,\frac{1}{2i}\{p,\bar{p}\}\}(\rho)=-2\times 3^{2} (|\im z|)^{4/3}\ne 0.\]
Nous sommes donc dans le cas 1.

Pour tout $z$ appartenant à l'intérieur de  $\Sigma(p),$ nous avons
\[ p^{-1}(z)=\{\rho_{+},\rho_{-}\},\]
où 
\[ \rho _{\pm}(z)=((\im z)^{\frac{1}{3}},
\mp(\re z)^{\frac{1}{2}}),\]
et 
\[\frac{1}{2i}\{p,\bar{p}\}(\rho_{\pm})
=(\im p_{\xi}\bar{p}_{x})(\rho_{\pm})=-2\times3 \xi x^{2}.\]
Le changement de variable $x=(\im z)^{\frac{1}{3}}\tilde{x},$ 
nous permet d'identifier
$(\mathcal{P}-z)$ à $(\im z)\widetilde{\mathcal{P}}$ où  $\widetilde{\mathcal{P}}$ s'écrit
\begin{equation}
\big(\frac{1}{(\im z)^{\frac{5}{6}}}D_{x}\big)^{2}
+ix^{3}-\frac{\re z}{\im z}-i,\mbox{ avec }
h:=\frac{1}{(\im z)^{\frac{5}{6}}},\;
y:=\frac{\re z}{\im z}.
\end{equation}
$h$ est notre paramètre semiclassique et $y$ est choisi petit. 
On choisit $y>0,$ puisque l'on se restreint aux courbes situées au dessus de l'axe réel.
Soit $\varphi(\xi)$ la fonction qui vérifie 
\[\xi^{2}+i\varphi(\xi)^{3}-y-i=0.\]
Nous avons
\[\int_{\sqrt{y}}^{-\sqrt{y}}\varphi d\xi=
\int_{\sqrt{y}}^{-\sqrt{y}}(1-iy+i\xi^{2})^{\frac{1}{3}}d \xi\]
Puisque $y$ est petit, un développement limité donne 
\begin{align*}
\im\int_{\sqrt{y}}^{-\sqrt{y}}(1-iy+i\xi^{2})^{\frac{1}{3}}d \xi=
\frac{4}{9} y^{\frac{3}{2}}
+\MO(y^{5/2}).
\end{align*}
 Ensuite, remarquons que
 \begin{align*}
 \frac{1}{2i}\{p,\bar{p}\}(\rho_{+}(y+i))&=
 \frac{1}{2i}\{p,\bar{p}\}(1,-y^{1/2})\\
 &=(-3x^{2}\xi-3x^{2}\xi)
 (1,-y^{1/2})
 =6\sqrt{y},
 \end{align*}
il en découle que
 \[ \frac{1}{2i}\{\bar{p},p\}(\rho_{-}(y+i))=6\sqrt{y}.\]
 \begin{prop}\label{cub1}
 Pour $\im z\gg 1$ et $y:=\frac{\re z}{\im z},$ $h:= \frac{1}{(\im z)^{\frac{5}{6}}},$ tels que
 \[ h\ll y^{3/2}\ll \MO(h\ln\frac{1}{h}) \mbox{ soit, } (\im z)^{4/9}\ll \re z\ll \im z, \]
 nous avons de manière précise,
 \[\|(\mathcal{P}-z)^{-1}\|\sim \frac{\sqrt{\pi}
 \exp\left(\frac{1}{h}\im \int_{\sqrt{y}}^{-\sqrt{y}}(1-iy+ix^{2})^{\frac{1}{3}}
 dx \right)}{
 6^{1/2}(\im z)^{1/3}(\re z)^{1/4}}
 +\MO(\frac{1}{(\re z)^{\frac{1}{4}}(\im z)^{\frac{1}{3}}}),\]
 ou bien, sous une forme simplifiée
\begin{equation}\label{cub.2} 
\|(\mathcal{P}-z)^{-1}\|\sim
\frac{\sqrt{\pi}\exp(\frac{4}{9}\frac{(\re z)^{3/2}}{(\im z)^{2/3}} 
+\MO(\frac{(\re z)^{5/2}}{(\im z)^{5/3}}))}
{ 6^{1/2}(\im z)^{1/3}(\re z)^{1/4}}
+\MO(\frac{1}{(\re z)^{\frac{1}{4}}(\im z)^{\frac{1}{3}}}).
\end{equation}
 \end{prop}
 \paragraph{Lignes de niveaux de la résolvante.}
On se sert de l'expression (\ref{cub.2}),
\[ \frac{1}{\epsilon}=\|(\mathcal{P}-z)^{-1}\|=(\ref{cub.2}).\]
On obtient alors 
\begin{equation}\label{cub.3}
-\ln \epsilon=\frac{4}{9}\frac{(\re z)^{3/2}}{(\im z)^{2/3}}+
\MO(\frac{(\re z)^{5/2}}{(\im z)^{5/3}})
+\ln\left( \frac{\sqrt{\pi}}{6^{1/2}}\frac{1}{(\re z)^{1/4}(\im z)^{1/3}}\right).
\end{equation}
(\ref{cub.3}) se met sous la forme
\begin{align}
\re z=&(\frac{9}{4})^{2/3} (\im z)^{4/9}
\left(\ln (\frac{6^{1/2}(\re z)^{1/4}(\im z)^{1/3}}
{\sqrt{\pi}\epsilon})\right)^{2/3}
+\MO(\frac{(\re z)^{5/3}}{(\im z)^{2/3}})\label{cub.4}
\nonumber
\\
=&(\frac{9}{16})^{2/3} (\im z)^{4/9}
\left(\ln (\re z( \im z)^{4/3})\right)^{2/3}
+\MO((\im z)^{4/9})
\nonumber\\
&+\MO(\frac{(\re z)^{5/3}}{(\im z)^{2/3}}).
\end{align}
Donc
\begin{align}
\re z=&(\frac{9}{16})^{2/3}(\im z)^{4/9}\left(\ln (\im z)^{16/9}
(\ln ((\re z)(\im z)^{4/3})^{2/3} 
\right)^{2/3}\nonumber \\
&+\MO((\im z) ^{4/9}).
\end{align}
Puisque $\re z<\im z$ et que $\ln(t\ln t)=(1+o(1))\ln t,$
on obtient
\begin{align}
\re z=(1+o(1))(\frac{9\cdot16}{16\cdot9})^{2/3}(\im z)^{4/9}(\ln \im z)^{2/3}.
\end{align}
Résumons :
\begin{prop}
La partie des lignes de niveaux de l'opérateur non-autoadjoint 
$\mathcal{P}$ se situant
en dessous de l'axe réel, admet la 
représentation asymptotique suivante lorsque $\im z\to \infty$
\begin{align*}
\re z=(1+o(1))
(\im z)^{4/9}(\ln \im z)^{2/3}.
\end{align*}
De manière plus précise, si l'on considère la ligne de niveaux 
$\frac{1}{\epsilon}=\|(\mathcal{P}-z)^{-1}\|$ alors nous avons 
\[
\re z=(1+o(1))
(\frac{9}{4})^{2/3} (\im z)^{4/9}
\left(\ln (\frac{(\im z)^{4/9}}{\epsilon})\right)^{2/3}.
\]
\end{prop}
  
\subsection{Cas 2: Oscillateur harmonique non-autoadjoint}
Maintenant, considérons l'opérateur harmonique non-autoadjoint,
\begin{equation}
Q=(D_{x})^{2}+ix^{2}
\end{equation}
sur la droite réelle.

$\Sigma(p)$ est le premier quadrant supérieur, et $\partial\Sigma(p)=
\R_{+}\cup i\R_{+}.$

Les lignes de niveaux de la résolvante sont symétriques 
par rapport à l'axe $e^{i\pi/4}\R.$ Intéressons-nous ici au comportement 
asymptotique des lignes de niveaux. Au regard de la symétrie axiale, 
il suffit de considérer la partie des lignes qui se situe en dessous 
de l'axe de symétrie. Nous avons le résultat suivant de K.~Pravda~Starov:
\begin{prop}
$S_{i}$ est la symétrie axiale par rapport à l'axe $e^{i\pi/4}\R.$  Il existe 
$C>0$ tel que nous avons
\begin{multline*}
\{z\in \C;\;\re z\ge \nu_{0}\mbox{ et }
0\le \im z<C(\re z)^{1/3}-\epsilon\}\cup\\
S_{i}\left(
\{z\in \C;\;\re z\ge \nu_{0}\mbox{ et }
0\le \im z<C(\re z)^{1/3}-\epsilon\}
\right) \subset \sigma_{\epsilon}(Q)^{c},
\end{multline*}
$\sigma_{\epsilon}$ désigne le $\epsilon$-pseudospectre de $Q.$
\end{prop}

Pour tout  $z$ dans l'image du symbole $q$, nous avons 
\begin{align}
q^{-1}&(z)=\{\rho_{-}^{j}(z),\rho_{+}^{j}(z)|\;j=1,2\},\\
&\mbox{où } \frac{1}{2i}\{q,\bar{q}\}(\rho_{\pm})>0.
\end{align}
et si $z_{0}$ un point du bord, appartenant à $\R_{+}$
alors
\[q^{-1}(z_{0})=\{(\xi_{0},0)\}\cap\{(-\xi_{0},0)\},\quad\xi_{0}>0.\]
Nous nous trouvons donc dans le cas 2.
Dans un voisinage de $z_{0}$ nous avons 
\begin{align*}
\rho_{+}^{1}(z)=((\im z)^{\frac{1}{2}},-(\re z)^{\frac{1}{2}}),\quad
&\rho_{-}^{1}(z)=(- (\im z)^{\frac{1}{2}},-(\re z)^{\frac{1}{2}}),\\
\rho_{+}^{2}(z)=(-(\im z)^{\frac{1}{2}},(\re z)^{\frac{1}{2}}),\quad
&\rho_{-}^{2}(z)=((\im z )^{\frac{1}{2}},(\re z)^{\frac{1}{2}}).
\end{align*}
Le changement de variable $x=(\re z)^{\frac{1}{2}}\tilde{x},$ nous permet d'identifier
$Q-z$ à $(\re z)\widetilde{Q}$ où $\widetilde{Q}$ s'écrit
\begin{equation}
\left(\frac{1}{\re z}D_{x}\right)^{2}+ix^{2}-1-i\frac{\im z}{\re z},\mbox{ avec }
h:=\frac{1}{\re z},\;
y:=\frac{\im z}{\re z}.
\end{equation}
$h$ est notre paramètre semiclassique et $y$ est choisi petit.

Nous avons 
\begin{align}
\xi^{2}+ix^{2}-1-iy=&(\xi-\sqrt{1+i y-ix^{2}})(\xi+\sqrt{1+i y-ix^{2}}),\nonumber\\
=:&(\xi-\varphi_{1}(x,y))(\xi-\varphi_{2}(x,y)).
\end{align}
Nous obtenons
\begin{align}\label{ex.6}
\int^{\sqrt{y}}_{-\sqrt{y}} \varphi_{1}\,dx=
\int^{-\sqrt{y}}_{\sqrt{y}} \varphi_{2}\,dx=
\int^{\sqrt{y}}_{-\sqrt{y}} \sqrt{1+i y-ix^{2}} dx.
\end{align}
Nous souhaitons calculer cette intégrale.

Soit $F(x)=\frac{1}{2} \arcsin x+\frac{x\sqrt{1-x^{2}}}{2},$ $F'(x)=\sqrt{1-x^{2}},$
nous avons
\begin{align}
 \frac{d}{dx}F\left(\frac{e^{i\frac{\pi}{4}}x}{\sqrt{1+iy}}\right)=
 \sqrt{1-\frac{ix^{2}}{1+iy}}\,\frac{e^{i\frac{\pi}{4}}}{\sqrt{1+iy}}.
 \end{align}
 Donc
 \[e^{-\frac{i\pi}{4}}(1+iy) \frac{d}{dx}F\left(\frac{e^{i\frac{\pi}{4}}x}{\sqrt{1+iy}}\right)
 =\sqrt{1+i y-ix^{2}}.\]
 Si bien que 
 \begin{align}
 \int^{-\sqrt{y}}_{\sqrt{y}} \varphi_{2}\,dx&=
 e^{-\frac{i\pi}{4}}(1+iy)\left(F\left(\frac{e^{i\frac{\pi}{4}}\sqrt{y}}{\sqrt{1+iy}}\right)
-F\left(\frac{-e^{i\frac{\pi}{4}}\sqrt{y}}{\sqrt{1+iy}}\right)\right)\nonumber\\
&= e^{-\frac{i\pi}{4}}(1+iy)\arcsin (\frac{e^{i\frac{\pi}{4}}\sqrt{y}}{\sqrt{1+iy}}),\label{ex.7}
\end{align}
un développement limité donne
\begin{align}
&=e^{-\frac{i\pi}{4}}(1+iy)(\frac{e^{i\frac{\pi}{4}}\sqrt{y}}{\sqrt{1+iy}}+
\frac{1}{6}\frac{e^{i\frac{3\pi}{4}}y^{3/2}}{(1+iy)^{3/2}}+\ldots),\\
&= \sqrt{1+iy}\sqrt{y} +\frac{1}{6} e^{i\frac{\pi}{2}}\frac{1}{(1+iy)^{\frac{1}{2}}}y^{3/2}+
\MO(y^{5/2}) .\label{ex.8}
\end{align}
\begin{remarque}
On aurait pu aussi traiter l'intégrale (\ref{ex.6}) en remarquant que 
\[\int\sqrt{x^{2}+a}=\frac{1}{2}(x\sqrt{x^{2}+a}+a\ln (x+\sqrt{x^{2}+a})),\]
ce qui donne 
\[\int^{-\sqrt{y}}_{\sqrt{y}} \varphi_{2}\,dx=
\frac{1}{2}e^{-\frac{\pi i }{4}}(i-y)\ln 
\frac{\sqrt{y}+\sqrt{i}}{-\sqrt{y}+\sqrt{i}}.\]
On retrouve (\ref{ex.7}) par la formule :
\[\arcsin(z)=-i\ln \left(iz+(1-z^{2})^{\frac{1}{2}}\right).\]
\end{remarque}
Finalement, en prenant la partie imaginaire dans (\ref{ex.8}), on obtient  
\[\im\int^{\sqrt{y}}_{-\sqrt{y}} \varphi_{1}\,dx=
\im \int^{-\sqrt{y}}_{\sqrt{y}} \varphi_{2}\,dx= \frac{1}{2}y^{3/2}
+\frac{1}{6}y^{3/2}=\frac{2}{3}y^{3/2}+\MO(y^{5/2}).\]
Ensuite, remarquons 
\begin{align*}
\frac{1}{2i}\{q,\bar{q}\}(\rho_{+}^{1}(1+iy))
&=\frac{1}{2i}\{q,\bar{q}\}(\sqrt{y},-1)\\
&=(-2\xi x-2\xi x)(-1,\sqrt{y})=4\sqrt{y},\\
\frac{1}{2i}\{\bar{q},q\}(\rho_{-}^{1}(1+iy))&=4\sqrt{y}.
\end{align*}

On déduit  alors le résultat suivant:
\begin{prop}
 Pour $\re z\gg 1$ et $y:=\frac{\im z}{\re z},$  $h:=\frac{1}{\re z},$
\[h\ll y^{3/2} \ll 1\mbox{ soit } (\re z)^{1/3}\ll \im  z\ll \re z,\] 
nous avons de manière précise,
\begin{align*}
\|(Q-z)^{-1}\|\sim 
\frac{\sqrt{\pi}
\exp\left(\frac{1}{h} \im (e^{-\frac{i\pi}{4}}(1+iy)\arcsin (\frac{e^{i\frac{\pi}{4}}\sqrt{y}}{\sqrt{1+iy}}))\right)}
{2(\re z)^{1/4}(\im z) ^{1/4}}\\
+\MO(\frac{1}{(\im z)^{\frac{1}{4}}(\re z)^{\frac{1}{4}}}),
\end{align*}
ou bien dans une forme symplifiée,
\begin{equation} \label{ex.12}
\|(Q-z)^{-1}\|\sim 
\frac{\sqrt{\pi}\exp(\frac{2}{3}\frac{(\im z)^{3/2}}{(\re z)^{1/2}}
+\MO(\frac{(\im z)^{5/2}}{(\re z)^{3/2}}))}
{2(\re z)^{1/4}(\im z) ^{1/4}}+\MO(\frac{1}{(\im z)^{1/4}(\re z)^{1/4}}).
\end{equation}
\end{prop}
\paragraph{Lignes de niveaux de la résolvante.}
On se sert de l'expression (\ref{ex.12}),
\[ \frac{1}{\epsilon}=\|(Q-z)^{-1}\|=(\ref{ex.12}).\]
On obtient alors 
\begin{equation}\label{ex.13}
-\ln \epsilon=\frac{2}{3}\frac{(\im z)^{3/2}}{(\re z)^{1/2}}+
\MO(\frac{(\im z)^{5/2}}{(\re z)^{3/2}})
+\ln\left( \frac{\sqrt{\pi}}{2}\frac{1}{(\re z)^{1/4}(\im z)^{1/4}}\right).
\end{equation}
En oubliant le facteur logarithmique dans (\ref{ex.13}), on aurait
\begin{align}\label{ex.14}
-\ln \epsilon&=\frac{2}{3}\frac{(\im z)^{3/2}}{(\re z)^{1/2}}+
\MO(\frac{(\im z)^{5/2}}{(\re z)^{3/2}})=
&\frac{2}{3}\frac{(\im z)^{3/2}}{(\re z)^{1/2}}(1+
\MO(\frac{\im z}{\re z}))
\end{align}
Soit 
\[(-\ln \epsilon)^{2}\frac{9}{4}=\frac{(\im z)^{3}}{\re z}(1+
\MO(\frac{\im z}{\re z})),\]
ainsi
\begin{equation}
\im z=(\re z\,(-\ln \epsilon)^{2}\,\frac{9}{4})^{1/3}(1+
\MO(\frac{\im z}{\re z})).
\end{equation}
(\ref{ex.13}) se met sous la forme
\begin{align}
\im z&=(\frac{3}{2})^{2/3} (\re z)^{1/3}
\left(\ln (\frac{2(\re z)^{1/4}(\im z)^{1/4}}{\sqrt{\pi}\epsilon})\right)^{2/3}
+\MO(\frac{(\im z)^{5/3}}{(\re z)^{2/3}})
\nonumber\\
&=(\frac{1}{6})^{2/3} (\re z)^{1/3}
\left(\ln (\re z\im z)\right)^{2/3}
+\MO((\re z)^{1/3})
+\MO(\frac{(\im z)^{5/3}}{(\re z)^{2/3}}).
\label{ex.15}
\end{align}
Donc
\begin{align}
\im z=&(\frac{1}{6})^{2/3}(\re z)^{1/3}\left(\ln (\re z)^{4/3}(\ln (\re z\im z)^{2/3} 
\right)^{2/3}
\nonumber\\
&+\MO((\re z)^{1/3}).
\end{align}
Puisque $\im z<\re z$ et que $\ln(t\ln t)=(1+o(1))\ln t$
on obtient
\begin{align}
\im z=(1+o(1))(\frac{1\cdot4}{6\cdot3})^{2/3}(\re z)^{1/3}(\ln \re z)^{2/3}.
\end{align}
Résumons :
\begin{prop}
La partie des lignes de niveaux de l'oscillateur harmonique non-autoadjoint se situant
en dessous de la droite $e^{i\pi/2}\R$ 
admet la 
représentation asymptotique suivante lorsque $\re z\to \infty$
\begin{align*}
\im z=(1+o(1))(\frac{1\cdot4}{6\cdot3})^{2/3}(\re z)^{1/3}(\ln \re z)^{2/3}.
\end{align*}
De manière plus précise si l'on considère la ligne de niveau 
$\frac{1}{\epsilon}=\|(Q-z)^{-1}\|$ alors nous avons 
\[
\im z=(1+o(1))
(\frac{3}{2})^{2/3} (\re z)^{1/3}
\left(\ln (\frac{(\re z)^{1/3}}{\epsilon})\right)^{2/3}.
\]
\end{prop}
On remarquera que les lignes de niveaux de la résolvante ne se 
rapprochent pas à l'infini puisque partant de (\ref{ex.15}), nous avons 
\begin{align*}
\im z=&(\frac{3}{2})^{2/3}(\re z)^{1/3}\left(\ln (\re z)^{4/3}(\ln (\re z\im z)^{2/3} 
-\ln \epsilon \right)^{2/3}\\
&+\MO(\frac{(\im z)^{5/3}}{(\re z)^{2/3}}), \\
=&C_{0}(\re z)^{1/3}(\ln \re z)^{2/3}\left(1+C_{1}\frac{-\ln \epsilon}{\ln \re z}+
\MO(\frac{\ln \epsilon}{\ln \re z})^{2}\right)\\
&+\MO(\frac{(\im z)^{5/3}}{(\re z)^{2/3}}), \\
=&C_{0}(\re z)^{1/3}(\ln \re z)^{2/3}+C_{0}'(\re z)^{1/3}\frac{-\ln \epsilon}{(\ln \re z)^{1/3}}+\ldots
\end{align*}
Le second membre tend vers l'infini pour un $\epsilon$ donné.\\

\noindent\textbf{Projection spectrale.}
Ici nous nous intéressons au projection spectrale 
de l'oscillateur harmonique non autoadjoint. 
\begin{theo}[Davies, Kuijlaars]
Soit l'opérateur $H=D_{x}^{2}+z^{4}x^{2}$ pour un complexe $z=e^{i\theta}$.
Les valeurs propres sont notées $\lambda_{n}=z^{2}(2n+1).$ On note 
$N_{n,z}$ la norme du projecteur spectral $\frac{1}{2\pi i} \int_{D(\lambda_{n},\epsilon)}(z-H)^{-1}dz.$
Si $0<\theta<\pi/4$ alors
\[\lim_{n\to \infty}n^{-1}\ln N_{n,z}=2\re (f(r(\theta)e^{i\theta}))\]
où
\[f(z)=z\sqrt{z^{2}-1}+\ln (z+\sqrt{z^{2}-1}),\]
et 
\[r(\theta)=(2\cos(2\theta))^{-1/2}.\]
\end{theo}

Nous allons redémontrer ce théorème avec nos estimations.\\
\noindent\textbf{Preuve.}
Le changement de variable $x=(2n+1)^{1/2}\tilde{x}$ permet d'identifier
$H-\lambda_{n}$ à 
$(2n+1)((hD_{x})^{2}+z^{4}x^{2}-z^{2}),$ avec $h=(2n+1)^{-1}.$
Compte tenu de qui précède, nous avons  
\begin{align*}
\lim_{n\to \infty}\ln N_{n,z}&=-\frac{1}{h}\im
\int^{x_{+}}_{x_{-}} \sqrt{(z^{2}-z^{4}x^{2})}\;dx
=-\frac{1}{h}\im
\int^{x_{+}}_{x_{-}} \sqrt{(1-z^{2}x^{2})}\;dx\\
&=\frac{1}{h}\re
\int^{x_{+}}_{x_{-}} \sqrt{(z^{2}x^{2}-1)}\;dx,
\end{align*}
ou
\begin{align*}
 x_{+}(z)&=\frac{\im z^{2}}{\im z^{4}}=\sqrt{\frac{\sin (2 \theta)}{\sin (4\theta)}}=\\
 &=\sqrt{\frac{\sin (2 \theta)}{2\sin(2\theta) \cos( 2 \theta)})}=(2\cos(2\theta))^{-1/2},\\
  x_{-}(z)&=-(2\cos(2\theta))^{-1/2}.
  \end{align*}
 Etant donné la remarque 4.4, il est facile de conclure.
 \hfill $\square$ \medskip

\subsection{Opérateur d'advection-diffusion}
Ici nous calculons la résolvante d'un 
opérateur d'advection-diffusion au point $1+i.$
 \\
On considère l'opérateur non-autoadjoint sur le cercle
\[L=-\sin(x)(hD_{x})^{2}-ihD_{x}.\]

Fort des informations de \cite{TE}, nous savons que l'image du symbole 
à l'infini est comprise entre les deux paraboles $\re z=\pm (\im z)^{2}.$
Nous allons montrer que c'est exactement 
la frontière de  son image. Pour cela, remarquons que 
\begin{align}
\ell^{-1}(z)&=\{(\arcsin(-\frac{\re z}{(\im z)^{2}}) ,-\im z)\}\cup
\{(\pi -\arcsin(-\frac{\re z}{(\im z)^{2}}) ,-\im z)\}\nonumber\\
&=\{\rho_{+}\} \cup \{\rho_{-}\}
\end{align}
Nous avons 
\begin{equation}\label{dd.1}
\frac{1}{2i}\{\ell,\bar{\ell}\}(\rho)=\xi^{2}\cos(x).
\end{equation}
Par souci de simplification, nous notons $\gamma:=\frac{\re z}{(\im z)^{2}}.$
Puisque $\cos(\arcsin x)=\sqrt{1-x^{2}}.$ nous voyons que 
\begin{align*}
\frac{1}{2i}\{\ell,\bar{\ell}\}(\rho_{+})&=(\im z)^{2}\sqrt{1-\gamma^{2}}>0,\\
\frac{1}{2i}\{\ell,\bar{\ell}\}(\rho_{-})&=-(\im z)^{2}\sqrt{1-\gamma^{2}}<0.
\end{align*}
En utilisant l'annulation du crochet de Poisson sur le bord de l'image de $l$ 
nous obtenons:
\begin{lemme}
Soit $L=-\sin(x)(hD_{x})^{2}-ihD_{x}$ de symbole $\ell$ alors 
$\ell(T^{\ast}{S^{1}})$ est bornée par les paraboles
\[\re z=-(\im z)^{2} \mbox{ et } \re z=+(\im z)^{2}.\]
De plus, nous affirmons que les points du bord de $\Sigma(\ell)$ sont tous d'ordre 2, à l'exception 
de l'origine.
\end{lemme}
\noindent\textbf{Preuve.}
Nous allons calculer le deuxième crochet pour vérifier que les point du bord sont d'ordre 2
\begin{align*}
\{\ell,\frac{1}{2i}\{\ell,\bar{\ell}\}\}(\rho)=&
\{-\sin(x)\xi^{2}-i\xi,\xi^{2}\cos(x))\}\nonumber\\
=&2\xi^{3}(\sin x)^{2}+2\xi^{3}(\cos x)^{2}\\
&+i\xi^{2}\sin x\\
=&2\xi^{3}+i\xi^{2}(\sin x ).
\end{align*}
Dans le cas où $\im z$ ne s'annule pas,
il est clair que le bord est d'ordre 2.
\hfill $\square$ \medskip




Soit le point $(1+i)-\alpha e^{i\theta}$ $\theta\in]-\pi/2,\pi/2[$ où 
$\cos\theta =-1/\sqrt{5},$  $\sin\theta =2/\sqrt{5}.$ 
Nous partons pour $\alpha$ petit dans la direction donnée par la
normale à la courbe $\partial\Sigma(\ell)$ au point $1+i.$

Nous avons les formules 
exactes plus les développements limités suivants (en utilisant 
$\arcsin(1-x^{2})\sim\frac{\pi}{2}+\sqrt{2}|x|+\MO(|x|^{3}))$:  
\begin{align*}
x_{+}(\alpha)&=\arcsin(-\frac{1-\alpha \cos \theta}{(1-\alpha \sin \theta)^{2}}),\\
&=-\frac{\pi}{2}+\alpha^{1/2}(\sqrt{2}\times(-\cos\theta+2\sin\theta))^{1/2}+\MO(\alpha)\\
&=-\frac{\pi}{2}+\alpha^{1/2}(\sqrt{2\times 5})^{1/2}+\MO(\alpha),\\
x_{-}(\alpha)&=\pi-\arcsin(-\frac{1-\alpha \cos \theta}{(1-\alpha \sin \theta)^{2}})\\
&=-\frac{\pi}{2}-\alpha^{1/2}(\sqrt{2\times 5})^{1/2}+\MO(\alpha).\\
\end{align*}
En utilisant la formule \ref{gg.6} on aurait retrouvé
\[x_{\pm}=-\frac{\pi}{2}\pm\sqrt{5}\alpha^{1/2}(\frac{\sqrt{2}}{\sqrt{ 5}})^{1/2}+\MO(\alpha).\]
Pour $\xi_{\pm}$ nous avons 
\[\xi_{\pm}(\alpha)=-1+\alpha \sin\theta.\]
Avec les données exactes de $x_{\pm}$ et $\xi_{\pm}$ nous retrouvons
\[\frac{1}{2i}\{\ell,\bar{\ell}\}(\rho_{\pm})\asymp \alpha^{1/2}.\]

Pour l'intégrale d'action, 
la fonction implicite $\varphi(x,z),$ vérifiant $\varphi(-\pi/2,1+i)=-1,$
de l'équation  
$-(\sin x) \xi^{2}-i\xi-z=0$ est
\begin{equation}
\frac{1}{2}\times \frac{-i+(-1-4 z \sin(x))^{1/2}}{\sin x},
\end{equation}
et 
\begin{align*}
\varphi'_{x}(x,z)=&
\frac{-i\cos x}{2(\sin x)^{2}}\\&+\frac{2z(\cos x)(-1-4z\sin x)^{-1/2}-
(\cos x)(-1-4z\sin x)^{1/2}}{(\sin x)^{2}}.
\end{align*}
Comme $\varphi'_{x}(-\pi/2,z)=0$, et 
$\varphi''_{x}(-\pi/2,z)\ne 0$ on peut retrouver 
\begin{equation}
\ell_{0}(z)=\im \int_{x_{-}}^{x_{+}}\varphi(x,z)dx\asymp \alpha^{3/2}.\\
\end{equation}

Nous avons alors redémontré dans le cas de l'opérateur d'advection-diffusion, 
la formule que nous connaissions \cite{SJ}:
\begin{prop} Soit $\theta$ tel que $\cos\theta =-1/\sqrt{5},$  $\sin\theta =2/\sqrt{5}.$
Pour  $h\ll\alpha^{3/2}<1,$ nous avons 
\begin{align*}
\|L-((1+i)-\alpha e^{i\theta})\|^{-1}\sim
\frac{1}{\sqrt{h}}
&\frac{\exp(\frac{1}{h}[\im\int\varphi(x,(1+i)-\alpha e^{i\theta}))]^{x_{+}}_{x_{-}})}{(1-\alpha \sin \theta )\sqrt{
1-\frac{1-\alpha \cos \theta}{(1-\alpha \sin \theta)^{2}}}}
(1+\MO(\frac{h}{\alpha^{3/2}}))\\
&+\MO(\frac{1}{\sqrt{h} \alpha^{1/4}}),
\end{align*}
où
\begin{equation*}
\varphi(x,z)=\frac{1}{2}\times \frac{-i+(-1-4 z \sin(x))^{1/2}}{\sin x},
\end{equation*}
avec
\begin{align*}
x_{+}(\alpha)&=\arcsin(-\frac{1-\alpha \cos \theta}{(1-\alpha \sin \theta)^{2}}),\\
x_{-}(\alpha)&=\pi-\arcsin(-\frac{1-\alpha \cos \theta}{(1-\alpha \sin \theta)^{2}}).\\
\end{align*}
Ce qui donne, après simplification
\[
\|L-((1+i)-\alpha e^{i\theta})\|^{-1}\asymp
\frac{C}{h^{1/2}\alpha^{1/4}}e^{\frac{\alpha^{3/2}}{Ch}}.\]
\end{prop}

\end {document}